\documentclass[12pt]{article}

\usepackage{amsmath}
\usepackage{amsfonts}
\usepackage{amsthm}

\newfont{\euler}{eufm10 scaled 1200} 

\def\Cmpx{{\mathbb{C}}}
\def\Real{{\mathbb{R}}}
\def\Intg{{\mathbb{Z}}}
\def\Rationals{{\mathbb{Q}}}
\def\natnum{{\mathbb{N}}}
\def\cnj{\overline}
\def\union{{\cup}}

\def\tfrac#1#2{{\textstyle{\frac{#1}{#2}}}}

\def\scrhalf{{\raisebox{.3ex}%
{$\scriptstyle 1$}\!/\!\raisebox{-.3ex}{$\scriptstyle 2$}}}

\def\innerprod#1{{\langle #1 \rangle}}
\def\norm#1{\|#1\|}
\def\rvec#1{{|#1\rangle}} 
 
\def\ad#1{{\textup{ad}_{#1}}}
\def\Wigner#1#2#3{{\left<
\begin{array}{ccc} & #2 & \\ #1 & & 0 \\ & #3 & \end{array}\right>}}
\def\CG#1#2#3#4#5#6{{C^{#1}_{#4}{}^{#2}_{#5}{}^{#3}_{#6}}}
\def\RMP#1#2#3{{{\cal R}^{#1}{}^{#2}{}^{#3}}}
\def\WsixJ#1#2#3#4#5#6{\left\{ \begin{array}{ccc} #1\ &\ #2\ &\ #3
\\#4\ &\ #5\ &\ #6\end{array}\right\}}

\def\defn#1{{\bf #1}}

\def\Set#1{\left\{#1\right\}}
\def\PB#1#2{\Set{#1,#2}}

\newtheorem{theorem}{Theorem}
\newtheorem{lemma}[theorem]{Lemma}

\newenvironment{jglist}
{\begin{list}{$\bullet$}
{
\setlength{\listparindent}{\parindent}\setlength{\parsep}{0 em}
}}
{\end{list}}

\newcounter{bean}
\newenvironment{countereg}
{\refstepcounter{bean}
\begin{trivlist}\item[]{\bf Counter example \arabic{bean}:}}
{\end{trivlist}}

\def\man{{\cal M}}

\def\Falg{{\cal F}}

\def\Falgfin{\Falg^{\scriptscriptstyle(\textup{fin})}}
\def\Falginf{\Falg^{{\scriptscriptstyle(}\infty{\scriptscriptstyle)}}}
\def\OFalg{{\cal F}^0}
\def\COFalg{{\cal F}^{C0}}
\def\Alg{{\cal A}}
\def\Algfin{\Alg^{\scriptscriptstyle(\textup{fin})}}
\def\Alginf{\Alg^{{\scriptscriptstyle(}\infty{\scriptscriptstyle)}}}
\def\CAlg{\Alg^C}
\def\OAlg{\Alg^0}
\def\SAlg{\Alg^\star}
\def\SOAlg{\Alg^{\star0}}
\def\CoM{{C^\omega(\man)}}
\def\CoSS{{C^\omega(S^2)}}
\def\CoMn#1{{C^\omega(\man_{#1})}}
\def\CoRn{{C^\omega(\Real^n)}}
\def\Algp{{\cal A_{\rho}}}

\def\elF#1{{\text{\euler #1}}}
\def\QC{\elF{c}}
\def\QI{\elF{b}}
\def\QCP{\elF{c}'}
\def\QIZ{\elF{b}^0}
\def\elFx{\elF{x}}
\def\elFy{\elF{y}}
\def\elFu{\elF{u}}
\def\elFv{\elF{v}}
\def\elFw{\elF{w}}

\def\elA#1{{\boldsymbol{#1}}}
\def\elAx{\elA{x}}
\def\elAX{\elA{X}}
\def\elAy{\elA{y}}
\def\elAY{\elA{Y}}
\def\elAu{\elA{u}}
\def\elAv{\elA{v}}
\def\elAw{\elA{w}}
\def\elAp{\elA{p}}
\def\elAq{\elA{q}}
\def\elAz{\elA{z}}

\def\elO#1{{#1}}
\def\elOx{\elO{x}}
\def\elOX{\elO{X}}
\def\elOy{\elO{y}}
\def\elOu{\elO{u}}
\def\elOv{\elO{v}}
\def\elOw{\elO{w}}
\def\elOp{\elO{p}}
\def\elOq{\elO{q}}

\def\elS#1{{\dot{#1}}}
\def\elSx{\elS{x}}

\def\elSy{\elS{y}}
\def\elSu{\elS{u}}
\def\elSv{\elS{v}}
\def\elSw{\elS{w}}

\def\q{Q}
\def\qC{Q^{\scriptscriptstyle{C}}}
\def\qI{Q^{\scriptscriptstyle{I}}}
\def\qIO{Q^{\scriptscriptstyle{I0}}}
\def\qCO{Q^{\scriptscriptstyle{C0}}}
\def\piC{\pi^{\scriptscriptstyle{C}}}
\def\piF{\pi^{\scriptscriptstyle{F}}}
\def\piS{\dot{\pi}}

\def\Hil{{\cal V}}
\def\Ideal{{\cal I}}
\def\Heis{{\cal H}}
\def\Omegas{{\tilde\Omega}}
\def\finite{{\scriptscriptstyle{\textup{finite}}}}
\def\hi{{\text{hi}}}
\def\lo{{\text{lo}}}

\def\Trn#1{\textup{tr}_{#1}}
\def\Tr{\textup{tr}}
\def\TrO{\textup{tr}_0}
\def\TrA{\textup{tr}_{\Alg}}

\def\dimman{{D}}

\title{Algebraic Noncommutative Geometry}
\author{{\bf Jonathan Gratus}\thanks{email: jg@luna.ph.lancs.ac.uk}\\
\small Physics Department, Lancaster University, Lancaster LA1 4YB}

\begin{document}

\maketitle

\begin{abstract}
A noncommutative algebra $\Alg$, called an algebraic noncommutative
geometry, is defined, with a parameter $\varepsilon$ in the
centre. When $\varepsilon$ is set to zero, the commutative algebra
$\OAlg$ of algebraic functions on an algebraic manifold $\man$ is
obtained.  This $\OAlg$ is a subalgebra of $\CoM$, which is dense if
$\man$ is compact.  The generators of $\Alg$ define an immersion of
$\man$ into $\Real^n$, and $\man$ inherits a Poisson structure as the
limit of the commutator. Thus $\Alg$ is a quantisation of a Poisson
manifold.  If an ordering convention is prescribed for $\Alg$ then a
star product on $\man$ is obtained.  Homomorphism and isomorphisms
between noncommutative geometries are defined, and the map from $\Alg$
to the Heisenberg algebra is used both to give an analogue of a
coordinate chart, and to give $\Alg$ a quantum group structure.
Examples of algebraic noncommutative geometries are given, which
include $\Real^n$, $T^\star S^2$, $T^2$, $S^2$ and surfaces of
rotation.  A definition of a metric on $\man$ which can be extended to
noncommutative geometry is given and this is used in an application of
noncommutative geometry to the numerical analysis of surfaces.
\end{abstract}

\tableofcontents


\section{Introduction}
\label{ch_intro}

\begin{table}[t]
\hspace{-0.8 cm}
\begin{tabular}{|c|c|c|c|c|c|c|c|}
\hline
Manifold &
Section &
Ordering &
Star &
Represent- &
Trace &
Heisenberg &
Quantum  
\\
& &  & product & ation &
(analytic) &
Coordinates &
 Group 
\\
\hline
Euclidean flat &
\ref{ch_egHeis} &
Wick &
Vey &
$\infty$-dimen- &
no &
yes &
yes
\\
\cline{3-4}
space $\Real^n$ & & 
Normal &
differential & sional
& & &
\\
\hline
Phase space &
\ref{ch_egTS2} &
& &
$\infty$-dimen-  &
no &
yes &
no
\\
$T^\star S^2$ & & & & sional & & &
\\
\hline
Torus $T^2$ or &
\ref{ch_egT2} &
Normal &
differential &
Matrix &
yes &
yes &
yes
\\
\cline{3-4}
Manin Plane &
&
Central &
Vey &
& (no) & &
\\
\hline
Surface of &
\ref{ch_egSR} &
Normal &
non-diff &
Matrix & 
yes & 
yes &
yes
\\
\cline{3-4}
Rotation &
&
Central &
Vey &
& (yes) & &
\\
\hline
Sphere $S^2$ &
\ref{ch_egS2} &
Normal &
non-diff &
Matrix &
yes &
yes &
yes
\\
\cline{3-4}
& &
Central &
Vey &
& (yes) & &
\\
\cline{3-4}
& &
Wick-like &
non-diff &
& & &
\\
\hline
Plane $\Real^2$ or $\Cmpx$ &
\ref{ch_egC} &
& & 
no &
no &
no &
no 
\\
\hline
\end{tabular}
\caption{The examples given in this article.}
\label{tbl_egs}
\end{table}

Noncommutative geometry has been suggested as a method for the
quantisation of gravity \cite{Kerner1}, string theory\cite{DeWit1},
renormalisation and a contribution to the elusive M-theory.  There are
two main goals in noncommutative geometry:

\begin{jglist}
\item 
Given a manifold or variety $\man$  the space of analytic functions
$\CoM$ forms an infinite dimensional commutative algebra via pointwise
multiplication. We wish to find a noncommutative algebra
$\Alg$ which can be considered as the noncommutative analogue of
$\CoM$. 

There are many possible principles to guide us to a definition of
$\Alg$. Here, guided by quantum mechanics, we define an element
$\varepsilon$ in the centre of $\Alg$, which plays the r\^{o}le of
$\hbar$. Thus when we set $\varepsilon=0$ we obtain a new commutative
algebra $\OAlg$ which is a (dense) subalgebra of $\CoM$.

\item We wish to define the tools of differential geometry, such as
tangent spaces, differential forms, connections and curvature, in
terms of the elements of $\CoM$, and then find analogues of these
objects when $\CoM$ is replaced by $\Alg$, so that they regain there
original definition when $\varepsilon=0$.
\end{jglist}

Most of this paper is concerned with the first of these goals, and,
having defined $\Alg$, giving detailed examples. The final section
gives an application for this theory to the numerical analysis of
surfaces.

The intrinsic method of defining a manifold is in terms of coordinate
charts. However the method employed here is to assume that $\man$ is
immersed in the real Euclidean space $\Real^n$. This implies that if
$\dim(\man)=\dimman$ then there are $n-\dimman$ functions
$\Set{\QIZ_1,\ldots,\QIZ_{n-\dimman}}$ with 
$\QIZ_s\colon\Real^n\mapsto\Real$
such that 
\begin{align}
\man=\Set{\underline{x}\in\Real^n\,|\,
\QIZ_s(\underline{x})=0,\,\forall \QIZ_s }
\end{align}
If the coordinates of $\Real^n$ are given by
$(\elOx_1,\ldots,\elOx_n)$ then each $\elOx_i\in\CoM$. These are, in a
certain sense, privilege elements of $\CoM$, as they encode all the
information about $\man$. They are called immersion coordinates.

In this article we shall further assume that $\man$ is algebraic; that
is, each $\QIZ_s(\underline{x})$ is a polynomial (multinomial) in the
coordinates $(\elOx_1,\ldots,\elOx_n)$.  Likewise we only consider the
subalgebra $\OAlg\subset\CoM$ of polynomials in
$(\elOx_1,\ldots,\elOx_n)$.  By restricting ourselves to algebraic
manifolds, we avoid many problems associated with
convergence. However, we shall see by the list of examples that this
still enables us to study a large class of interesting manifolds.

Since $\OAlg$ is a commutative algebra, we have the commutation
equations:
\begin{align}
[\elOx_i,\elOx_j]=0
\label{intro_Ox_comm}
\end{align}
where the square bracket represents the commutator.  Together with the
immersion equations $\Set{\QIZ_s(\underline{x})=0}$, this gives a
total of $n(n-1)/2+n-\dimman$ equations, which completely specify
$\OAlg$.

The noncommutative algebra $\Alg$ is also specified in this way.  It
is generated by the immersion coordinates
$\Set{\elAx_1,\ldots,\elAx_n}$ and a parameter $\varepsilon$ in the
centre of $\Alg$. We use the bold font to specify that
$\elAx_i\in\Alg$ and hence do not commute.  We replace
(\ref{intro_Ox_comm}) with
\begin{align}
[\elAx_i,\elAx_j]=i\varepsilon\elA{c}'_{i j}
\label{intro_Ax_comm}
\end{align}
for some $\elA{c}'_{i j}\in\Alg$. Thus when we set $\varepsilon=0$ this
reduces to (\ref{intro_Ox_comm}).  To make all this mathematically
precise we define everything by quotienting the algebra $\Falg$, which
is the free noncommuting algebra generated by
$(\elFx_1,\ldots,\elFx_n,\varepsilon)$, by various ideals.  The Euler
font is used for elements of $\Falg$. The details of this are in
section \ref{ch_defANCG}.

We call the algebra $\Alg$ an algebraic noncommutative geometry,
abbreviated to ANCG, to distinguish it from the noncommutative
geometry of Connes \cite{Connes1,Connes2} and the matrix geometry of
Madore \cite{Madore_bk}.

The act of setting $\varepsilon=0$ is given by the quotient map
$\pi\colon\Alg\mapsto\OAlg$. This may be considered as {\it
classicisation} i.e. taking one from a quantum system to a classical
system. The first consequence of this definition is that $\man$
inherits a Poisson structure. This is given, in section \ref{ch_Poi},
by
\begin{align*}
\PB{\pi(\elAu)}{\pi(\elAv})&=\pi(\tfrac1{i\varepsilon}[\elAu,\elAv])
\end{align*}
where $\elAu,\elAv\in\Alg$. Thus we have a method for the quantisation
of immersed Poisson manifolds.

\vspace{1 em} 

In section \ref{ch_Ord} we consider orderings.  Since
$\elAu\elAv\ne\elAv\elAu$, then, in order to specify the quantum
analogue of a particular classical function, we must specify an
ordering convention.  This is given by a linear map
$\Omega:\OAlg\mapsto\Alg$, so that $\pi\circ\Omega=1_{\OAlg}$.  In
quantum mechanics one often only has to specify an ordering for the
Hamiltonian.  One of the advantages with our approach to
quantisation, is that we do not, a priori, assume an ordering
convention, and can therefore compare different ordering conventions
on the same algebra $\Alg$.  For example for the Heisenberg algebra
$\Heis_2$ where $[\elAp,\elAq]=i\varepsilon$, which underlies the
quantum mechanics of a free particle on a line, we often consider two
orderings: 
\begin{align*}
& \text{Wick ordering;} 
&
\Omega_W (\elOp^r \elOq^s)&=\text{ sum of symmetric permutations
of $\elAp^r\elAq^s$} \\
& \text{Normal ordering;}
&
\Omega_N(\elOp^r \elOq^s)&=\elAp^r \elAq^s
\end{align*}
where $\pi(\elAp)=\elOp$ and $\pi(\elAq)=\elOq$.  Once we fix a
particular ordering then the algebra $(\Alg,\Omega)$ is equivalent to
a star algebra \cite{Flato1}; that is, we can transform the product on
$\Alg$ to a star product given by
\begin{align}
\elOu \star \elOv &= \sum_{r=0}^{\infty} \varepsilon^r C_r(\elOu,\elOv)
\label{intro_star}
\end{align}
where $C_r:\OAlg\times\OAlg\mapsto\OAlg$.  If one chooses the Wick
ordering of $\Heis_2$ then we have a Vey product where
$C_r=(\tfrac{i}2{\cal P})^r/r!$ where ${\cal
P}(\elOu,\elOv)=\PB{\elOu}{\elOv}$. On the other hand, if one chooses
the normal ordering of $\Heis_2$ then one has a different differential
star product.

Of course star products can be given intrinsically, simply by specifying
the functions \break $C_r\colon\CoM\times\CoM\mapsto\CoM$.
We have the following pseudo equation:
\begin{align}
\text{Algebraic Noncommutative Geometry}
+
\text{Ordering}
&\approx
\text{Star Product}
+
\text{Immersion Coordinates}
\label{intro_star_apprx}
\end{align}
This is not strictly true though, due to our restriction to algebraic
functions.

\vspace{1 em}

Completely separate to the question of which ordering to impose, is
the question of whether representations of $\Alg$ exist. This is
investigated in section \ref{ch_Rep}.  If $\man$ is compact then there
may exists a sequence of matrix representation of $\Alg$. These are
maps $\varphi_N\colon \Alg\mapsto M_n(\Cmpx)$, such that
$\varphi_N(\varepsilon)=\varepsilon_N{\mathbf 1}_N$, with
$\varepsilon_N\in\Real$, $\varepsilon_N\to0$ as $N\to\infty$. The
algebra of matrices which are the image of $\varphi_N$ can be thought
of as a matrix geometry.  For compact symplectic manifolds the limit
of the trace can be written in terms of the integral over $\man$.
\begin{align*}
\lim_{N\to\infty}\tfrac1N\text{tr}(\varphi_N(f))
&=
\frac1{|\man|}\int_\man \pi(f) \omega^n
\end{align*}
where $\omega$ is the symplectic 2-form and $|\man|=\int_\man \omega^n$.

\vspace{1 em}

In section \ref{ch_Hom} we define the concept of homomorphisms and
isomorphisms between noncommutative geometries.  An important case is
when the codomain is the Heisenberg algebra, which is the
noncommutative analogue of Euclidean flat space. This mapping may then
be considered as the noncommutative analogue of a coordinate system
(section \ref{ch_Heis}). We give examples of such noncommutative
coordinate systems. This returns us to the original ideas of Dirac who
suggested that one could consider manifolds where the coordinates do
not commute. Since the Heisenberg algebra can be given the structure
of a quantum group we can use the coordinate homomorphisms to give the
quantum group structure to other noncommutative geometries (section
\ref{ch_QG}).

\vspace{1 em}

In section \ref{ch_eg}, we give a number of examples. These are the
cotangent bundle of the sphere, flat space and the two dimensional
manifolds of the plane, torus, sphere, and surfaces of rotation. A
list of properties is given is table \ref{tbl_egs}.  The cotangent
bundle, section \ref{ch_egTS2}, should be thought of as the
noncommutative analogue of a phase space.  Thus we demonstrate that
algebraic noncommutative geometry is indeed a method of
quantisation. This gives the underlying quantum algebra corresponding
to the non-relativistic quantisation of a free particle on a sphere.

\vspace{1 em}


As mentioned the second goal of noncommutative geometry is to write
down the objects studied in differential geometry, such as tangent
bundles, cotangent bundles, exterior algebras, metric tensors,
connections and curvature, in terms of elements of the algebra
$\CoM$ and then find analogues of these objects when $\CoM$ is
replaced by $\Alg$, so that they regain there original definition when
$\varepsilon=0$.

There are two key properties required of tangent vector
fields. Firstly that they should be derivatives, i.e. follow Leibniz
rule, and secondly that they should form a module over the algebra of
functions. (That is one can multiply a vector with a scalar to give
another vector.)  It turns out that for noncommutative geometry these
two properties are incompatible, and one must choose either to have
vectors which are derivatives, or vectors which form a module.

The standard method is to choose vectors which form derivatives
\cite{Madore_bk}; that is, $\xi:\Alg\mapsto\Alg$ such that
$\xi(\elAu\elAv)=\xi(\elAu)\elAv+\elAu\xi(\elAv)$ for all
$\elAu,\elAv\in\Alg$.  For many noncommutative geometries we can show
that this implies that all vector are inner, i.e.  there exists
$\elAw\in\Alg$ such that $\xi=\tfrac1\varepsilon\ad{\elAw}$ where
$\ad{\elAw}(\elAu)=[\elAw,\elAu]$.
Clearly if $\xi$ is inner then $\elAu\xi$ is not inner.

In \cite{JG-TS2,JG-TSR} the author gives an alternative method of
defining tangent vectors on the noncommutative sphere and surfaces of
rotation. These vectors do form a (one sided) module over the
noncommutative surface but are derivatives only in the commutative
limit; that is,
$\xi(\elAu\elAv)=\xi(\elAu)\elAv+\elAu\xi(\elAv)+O(\varepsilon)$.

In this article, section \ref{ch_Geo}, we circumvent the problem of
how to define a vector field by defining the objects in differential
geometry using only the elements of $\CoM$. This we do by writing the
metric on $\man$ as $g({d\elOu}^\sharp,{d\elOv}^\sharp)$, where $g$
is the metric inherited from the ambient Euclidean immersion space and
$\sharp:T^\star\man\mapsto T\man$ is the metric dual.  We show that
this expression can be written using only the Poisson structure, the
immersion coordinates, and the local coordinates on $\man$.

\vspace{1 em}

In section \ref{ch_Mod}, we outline a method for the numerical
analysis of surfaces embedded in $\Real^3$.  Let us assume we wish to
analysis a surface, $\man$, which is nearly spherical; that is, the
function $\Set{\elOx_1,\elOx_2,\elOx_3}$, when expanded in spherical
harmonics, converges quickly. It therefore makes sense to use this
property in any numerical analysis of $\man$, and to encode the
information about the problem in terms of spherical harmonics as
opposed to pointwise encoding.

Let us assume we wish to calculate simply $u=vw$ where
$u,v,w\in\CoM$. We can express these functions using spherical
harmonics as $u=\sum_{nm}u_{nm}\psi^n_m$ etc. From the Echart-Wigner
theorem we have
\begin{align*}
u_{nm}=\sum_{m_1,n_1,n_2}
v_{n_1m_1} w_{n_2,m_2}
C^{n_1}_{m_1}{}^{n_2}_{m-m_1}{}^{n}_{m}
C^{n_1}_{0}{}^{n_2}_{0}{}^{n}_{0}
\end{align*}
Now if we work numerically then we truncate this sum and loose all
modes $\psi^n_m$ for $n\ge N$, for some $N\in\natnum$. As a result, in
general $(uv)w\ne u(vw)$. Thus the corresponding algebra, although it
is commutative, is non-associative.

By contrast, we propose, in section \ref{ch_Mod}, to use
the noncommutative spherical harmonics $\elA{P}^m_n$ described in
section \ref{ch_egS2}.  The corresponding algebra is
associative but noncommutative, indeed for numerical work this algebra
is simply the algebra of $N\times N$ matrices.  The method uses the
results described in this article so that differentiation is replaced
by commutation, and integration is replaced by trace.

All the information lost (or error) is introduced when we convert
functions $u$ on $\man$ into $N\times N$ matrices.  After this, we can
multiply any number of matrices without loosing additional
information.  Although the answer depends on the ordering of the
expression we wish to calculate, we will show that any difference will
be of order $O(1/N)$.


\vspace{1 em}

Finally in section \ref{ch_Disc} we discuss some of the possible
methods of enlarging $\Alg$ so that $\OAlg$ includes all analytic
function on $\man$. We also discuss other developments of this theory
and possible applications in physics.


\vspace{1 em}

{\bf NOTE}: This article is arranged so that all the theorems are
stated and proved before the main examples are given. This may not be
the easiest way to read this article and the casual reader is
recommended to scan the examples is section \ref{ch_eg} before and
whilst reading the theorems in section \ref{ch_ANCG}.


\subsection{Notation}

In this article we have a number of algebras, with many maps between
them. The elements in the main algebras are written with different
scripts to aid understanding. These are given by
\begin{align*}
\begin{array}{cccl}
\text{algebra} & \text{generators} & \text{general element} &
\\
\Falg & \Set{\varepsilon,\elFx_1,\ldots,\elFx_n} & 
\elFu,\elFv,\elFw,\elFy &
\\
\Alg & \Set{\varepsilon,\elAx_1,\ldots,\elAx_n} & 
\elAu,\elAv,\elAw,\elAy &
\\
\OAlg & \Set{\elOx_1,\ldots,\elOx_n} & 
\elOu,\elOv,\elOw,\elOy &
\\
\SAlg & \Set{\varepsilon,\elSx_1,\ldots,\elSx_n} & 
\elSu,\elSv,\elSw,\elSy &
\end{array}
\end{align*}

The expression $\CoM$ refers to the algebra of complex valued analytic
functions on $\man$. This means that for each function in $\CoM$ there
is a Taylor expansion about each point in $\man$. 

The term \defn{polynomial} in $(\elAx_1,\ldots,\elAx_n)$ means any
expression generated by taking finite sums and products. 

The square brackets always refer to the commutator, so that
$[\elAu,\elAv]=\elAu\elAv-\elAv\elAu$.

When talking about elements of $\Alg$ we will use the notation
$\elAu=O(\varepsilon^r)$ for $r\in\natnum$ to mean
$\elAu=\varepsilon^r \elAu'$ for some $\elAu'\in\Alg$. Likewise for
elements of $\Falg$ and $\SAlg$.

Algebraic noncommutative geometry is abbreviated to ANCG, whilst
ordered algebraic noncommutative geometry is abbreviated to OANCG.


\section{Definition and Properties of Algebraic Noncommutative
Geometries}
\label{ch_ANCG}

\subsection{Definition of ANCG}
\label{ch_defANCG}

Let $\OFalg$ be the free associative noncommutative complex algebra
with a unit, finitely generated by $\Set{\elFx_1,\ldots,\elFx_n}$. Let
$\Falgfin$ and $\Falginf$ be the algebra with elements
\begin{align}
\Falgfin=\Set{\sum_{r=0}^\finite \varepsilon^r\elFu_r
\,\bigg|\,\elFu_r\in\OFalg}
\qquad
\Falginf=\Set{\sum_{r=0}^\infty \varepsilon^r\elFu_r
\,\bigg|\,\elFu_r\in\OFalg}
\label{defANCG_exp_F}
\end{align}
where $\varepsilon$ is in the centre of $\Falgfin$ and $\Falginf$.
The notation $\displaystyle\sum_{r=0}^\finite$ means a finite sum over
non negative $r$. Unless otherwise specified, we write $\Falg$ to mean
either $\Falgfin$ or $\Falginf$.

Since $\varepsilon$ is in the centre of $\Falg$ we
can define the map $\piF\colon
\Falg\mapsto\OFalg=\Falg/\Set{\varepsilon\sim 0}$.  Thus $\piF$ is the
equivalent to setting $\varepsilon=0$. If $\elFu$ is written as
(\ref{defANCG_exp_F}) then $\piF(\elFu)=\elFu_0$.

We also define the quotient algebra and quotient map
\begin{align}
\qCO\colon \OFalg\mapsto\COFalg=\OFalg/\Set{[\elFx_i,\elFx_j]\sim0}
\end{align}
So $\COFalg$ is the free commutative algebra generated by
$\Set{\elOx_1,\ldots,\elOx_n}$ where $\elOx_i=\qCO(\piF(\elFx_i))$. We
specify that $\elOx_i=\cnj{\elOx_i}$ so that $\COFalg$ is a subalgebra
of $\CoRn$.

We define an \defn{Algebraic Noncommutative geometry} $\Alg$ as a
quotient algebra of $\Falg$. When we need to be precise we will write
$\Algfin$ or $\Alginf$ if it is the quotient of $\Falgfin$ or
$\Falginf$ respectively.

\begin{jglist}

\item $\Alg$ is noncommutative and associative.

\item $\Alg$ is the quotient of the algebra $\Falg$, for some
$n\in\natnum$, via the ideal generated from quotient elements
\begin{align}
&\QC_{ij}\in\Falg\,,\qquad i,j=1,\ldots,n
\\
&\QI_{s}\in\Falg\,,\qquad s=1,\ldots,n-\dimman
\end{align}
for some $\dimman\in\natnum$, $\dimman\le n$.  The ideal is all
elements of the form
\begin{align*}
\sum_{i,j=1}^n\elFu_{ij} \QC_{ij} + 
\sum_{s=1}^{n-\dimman}\elFv_s\QI_{s}\,,
\qquad\elFu_{ij},\elFv_s\in\Falg
\end{align*}
Since $\Alg$ is associative, this ideal must be a two sided ideal.
The quotient map is written 
\begin{align}
\q\colon \Falg\mapsto\Alg\,\text{ with }
\q(\varepsilon)=\varepsilon\,,\
\q(\elFx_i)=\elAx_i
\end{align}

\item Since $\varepsilon$ is in the centre of $\Alg$ we can quotient
$\Alg$ by the ideal generated from $\varepsilon$. This is equivalent
to setting $\varepsilon=0$. Thus we define the map
\begin{align}
\pi\colon \Alg\mapsto\OAlg=\Alg\big/\Set{\varepsilon\sim0}
\qquad\text{with }
\pi(\varepsilon)=0\,,\
\pi(\elAx_i)=\elOx_i
\end{align}

\item The commutation quotient relations $\QC_{ij}$ obey
\begin{align}
\QC_{ij}=\elFx_i\elFx_j-\elFx_j\elFx_i-i\varepsilon\QCP_{ij}
\qquad\text{where }
\QCP_{ij}\in\Falg
\label{defANCG_Cij}
\end{align}
and where there is at least one $\QCP_{ij}$ such that 
$\pi(\q(\QCP_{ij}))\ne0$.

\item The immersion quotient relations $\QI_{s}$ obey
\begin{align}
\qCO(\piF(\QI_{s}))\ne0 \qquad\forall s=1,\ldots,n-\dimman
\end{align}

\item $\Alg$ be a conjugation algebra; that is, there
exists a conjugation $\dagger\colon \Alg\mapsto\Alg$, where
\begin{align}
(\elAu\elAv)^\dagger=\elAv^\dagger{}\elAu^\dagger\,,
&&
\elAu^{\dagger\dagger}=\elAu\,,
&&
\varepsilon^\dagger=\varepsilon\,,
&&
(\elAx_i)^\dagger=\elAx_i\,,
&&
\lambda^\dagger=\cnj{\lambda}\,,
\qquad\forall\,\elAu,\elAv\in\Alg\,,\,\lambda\in\Cmpx
\end{align}

\end{jglist}

There are several more algebras, all of which are quotients of $\Falg$
which are useful. These are defined as follows, with their
corresponding quotient maps.
\begin{equation}
\begin{array}{rl}
&\qC \colon \Falg\mapsto\CAlg=\Falg/\Set{\QC_{ij}\sim0}
\\
&\qI \colon \CAlg\mapsto\Alg=\Alg/\Set{\QI_{s}\sim0} 
\qquad
\text{so that $\q=\qI\circ\qC$}
\\
&\qCO \colon  \OFalg\mapsto\COFalg=\Falg/\Set{[\elFx_i,\elFx_j]\sim0}
\\
&\qIO\colon \OFalg\mapsto\OAlg=\OFalg/\Set{\QIZ_{s}\sim0}
\qquad\text{where }
\QIZ_{s}=\piF(\qCO(\QI_{s}))
\\
&\piF\colon \Falg\mapsto\OFalg=\Falg/\Set{\varepsilon\sim0}
\\
&\piC\colon \CAlg\mapsto\COFalg=\CAlg/\Set{\varepsilon\sim0}
\end{array}
\label{defANCG_qot}
\end{equation}
The algebras $\Falg,\OFalg,\COFalg$ depend only on $n$, whilst
$\CAlg,\Alg,\OAlg$ depend on $n$ and the quotient relations $\QC_{ij}$
and $\QI_{s}$. Since all the maps simply correspond to quotients they
are related via the following commutative diagram.
\begin{align}
\begin{array}{ccccc}
\Falg & 
\stackrel{\qC}{\longrightarrow} &
\CAlg &
\stackrel{\qI}{\longrightarrow} &
\Alg
\\
\downarrow^{\piF} &&
\downarrow^{\piC} &&
\downarrow^{\pi}
\\
\OFalg &
\stackrel{\qCO}{\longrightarrow} &
\COFalg &
\stackrel{\qIO}{\longrightarrow} &
\OAlg
\end{array}
\label{defANCG_com_diag}
\end{align}
Since $\COFalg\subset\CoRn$, we can write
$\QIZ_{s}\colon\Real^n\mapsto\Cmpx$. Let
\begin{align}
\man = \Set{\underline{x}\in\Real^n\,\big|\,
\QIZ_s(\underline{x})=0, s=1,\ldots,n-\dimman}
\label{defANCG_QIZ}
\end{align}
In general $\man$ is an algebraic verity. If there are no critical
points then $\man$ is a manifold. Thus $\OAlg\subset\CoM$ is the
commutative subalgebra of complex algebraic function on $\man$.

Since $\elOx_i=\cnj{\elOx_i}$ for $\elOx_i\in\COFalg$ then
$\elOx_i=\cnj{\elOx_i}$ for $\elOx_i\in\OAlg$. It is easy to show that
this implies the $\pi$ preserves conjugation; that is,
$\pi(\elAu^\dagger)=\cnj{\pi(\elAu)}$ for $\elAu\in\Alg$.

\begin{lemma}
If $\man$ is compact then $\OAlg$ is dense in the space $C(\man)$ of
continuous complex valued functions on $\man$ with the uniform norm
\begin{align}
\norm{f}=\sup_{x\in\man}|f(x)|
\label{defANCG_norm}
\end{align}
\end{lemma}
\begin{proof}
Follows from the Boltzano-Wiestrass theorem.
\end{proof}

\subsection{Poisson Structure} 
\label{ch_Poi}

\begin{theorem}
There exists a Poisson structure on $\OAlg$ given by
\begin{align}
\PB{\bullet}{\bullet}\colon \OAlg\times\OAlg\mapsto\OAlg \,;
\qquad
\PB{{\pi(\elAu)}}{{\pi(\elAv)}}=\pi(\tfrac1{i\varepsilon}[\elAu,\elAv])
\label{Poi_def_Poi}
\end{align}
This can be extended to a Poisson structure on $\man$ given by
$\PB{\bullet}{\bullet}\colon  C^\infty(\man)\times C^\infty(\man)\mapsto
C^\infty(\man)$.
\end{theorem}

\begin{proof}
Given $\elAu,\elAv\in\Alg$, then since $\pi$ is a homomorphisms,
$\pi([\elAu,\elAv])=[\pi(\elAu),\pi(\elAv)]=0$. Thus
$[\elAu,\elAv]=O(\varepsilon)$, so
$\tfrac1\varepsilon[\elAu,\elAv]\in\Alg$. Hence the Poisson bracket is
defined. Given $\elAu_1,\elAu_2\in\Alg$ such that
$\pi(\elAu_1)=\pi(\elAu_2)$ then
$\elAu_1-\elAu_2=O(\varepsilon)$. Hence
$[\elAu_1-\elAu_2,\elAv]=O(\varepsilon^2)$, so
$\PB{\elOu_1}{\elOv}=\PB{\elOu_2}{\elOv}$, and the Poisson bracket is
well defined. From (\ref{defANCG_Cij}) we know there exist
$\elAu,\elAv\in\Alg$ such that
$\pi(\tfrac1{i\varepsilon}[\elAu,\elAv])\ne0$, and hence the Poisson
bracket is non trivial.

The derivative property follows from expanding
$\pi(\tfrac1{i\varepsilon}[\elAu,\elAv\elAw])$. The Jacobi identity
follows from the Jacobi identity for commutators.

Since all Poisson structures my be written in terms of a bivector,
the Poisson structure can be extended to $C^\infty(\man)$.
\end{proof}

We say that $\Alg$ is \defn{symplectic}, if the Poisson structure on
$\man$ is symplectic.  


\subsection{Orderings and Star Products}
\label{ch_Ord}

An \defn{ordering} on a ANCG $(\Alg,\man)$ is a choice of injective
linear map $\Omega\colon \OAlg\mapsto\Alg$ such that
$\pi\circ\Omega=1_{\OAlg}$ and $\Omega(1)=1$. Furthermore, $\Omega$ is
a \defn{unitary ordering} if
$\Omega(\cnj{\elOu})=\Omega(\elOu)^\dagger$.

The following theorem proves that all ANCG possess at least one
ordering. In general for a given ANCG there will be an infinite
number of orderings $\Omega$. 

\begin{theorem}
\label{th_Ord_exist_Ord}
Given a ANCG $(\Alg,\man)$ there exists a non unique unitary ordering
$\Omega$.
\end{theorem}

\begin{proof}
Choose a sequence of self-conjugate elements in $\Falg$ which is a
basis of $\OFalg$ as an $\infty$-dimensional vector space. For example
\begin{align}
\elFx_1,\elFx_2,\ldots,\elFx_n,
\elFx_1\elFx_1,\elFx_1\elFx_2,\ldots,\elFx_1\elFx_n,
\elFx_2\elFx_1,\ldots,\elFx_2\elFx_n,\ldots,
\elFx_n\elFx_1,\ldots \elFx_n\elFx_n,\elFx_1\elFx_1\elFx_1\ldots,\ldots
\label{Ord_PrincNorm}
\end{align}
now remove any elements from the sequence that are permutation of
previous elements, or are in the span of preceding elements and the
ideal generated by the quotients elements $\QIZ_{s}$. This gives a
sequence of $\elFu_i\in\Falg$. Each element in $\OAlg$ can be uniquely
written as a sum $\elOv=\sum_{i=0}^\finite \lambda_i\pi(\q(\elFu_i))$.
Set $\Omega(\elAv)= \sum_{i=0}^\finite \lambda_i\q(\elFu_i)$.

To construct a unitary ordering set
$\Omega_U(\elOv)=\tfrac12(\Omega(\elOv)+\Omega(\cnj{\elOv})^\dagger)$.

$\Omega$ is far from unique. For example we can always set a new
ordering as $\Omega_1(\elOv)=\Omega(\elOv)+\varepsilon$.
\end{proof}

There are certain orderings that have been given names.  For example
the \defn{Wick ordering} and \defn{normal ordering}. However the exact
definition of these orderings depend on the algebra $\Alg$.

\vspace{1 em}

We define an \defn{Ordered Algebraic Noncommutative Geometry} (OANCG)
as an ANCG $\Alg$ together with a choice of orderings $\Omega$. We can
now give the relation between OANCG and star products. We define the
set $\SAlg$ as the set of all elements of the form
\begin{align}
\SAlg = \Set{\elSu=\sum_{r=0}^\infty \varepsilon^r \elOu_r \bigg|
\quad \elOu_r\in\OAlg}
\end{align}

\begin{theorem}
We can extend the ordering $\Omega$ to give the map
\begin{align}
&\Omegas \colon  \SAlg  \mapsto \Alg\,; 
\qquad
\Omegas\left(\sum_{r=0}^\infty \varepsilon^r \elOu_r\right)=
\sum_{r=0}^\infty \varepsilon^r \Omega(\elOu_r)
\end{align}
This map has an inverse given by
\begin{align}
&\Omegas^{-1}\colon  \Alg\mapsto \SAlg\,;
\qquad
\Omegas^{-1}\colon  \elAu \mapsto \sum_{r=0}^\infty \varepsilon^r 
C_r(\elAu)
\end{align}
where $C_r\colon \Alg\mapsto\OAlg$ is given by
\begin{align}
C_r(\elAu)=\pi\left(\varepsilon^{-n}
\left(\elAu-\sum_{m=0}^{r-1} \varepsilon^m 
\Omega(C_m(\elAu))\right)\right)
\label{Ord_def_Cr1}
\end{align}
This satisfies $\Omegas\circ\Omegas^{-1}=1_{\Alg}$ and
$\Omegas^{-1}\circ\Omegas=1_{\SAlg}$
\end{theorem}
\begin{proof}
Trivial.
\end{proof}


Given an OANCG we can define a \defn{star product} on $\OAlg$. 
This is given by 
\begin{align}
\star\colon \OAlg\times \OAlg \mapsto \SAlg\,;
\qquad
\elOu\star \elOv = \Omegas^{-1}\Big(\Omega(\elOu)\Omega(\elOv)\Big)
=\sum_{r=0}^\infty \varepsilon^r C_r(\elOu,\elOv)
\label{Ord_def_star}
\end{align}
where
\begin{align}
C_r(\elOu,\elOv)=C_r(\Omega(\elOu)\Omega(\elOv))
\end{align}
We note that $C_0(\elOu,\elOv)=\elOu\elOv$ and
$C_1(\elOu,\elOv)-C_1(\elOv,\elOu)=i\PB{\elOu}{\elOv}$.  We extend
(\ref{Ord_def_star}) to the star product
$\star:\SAlg\times\SAlg\mapsto\SAlg$. We call the set $\SAlg$ together
with the product $\star$, a \defn{star algebra}.  This makes
$\Omegas:\Alg\mapsto\SAlg$ a bijective homomorphism. A
\defn{differentiable star product} requires that $C_r(\elOu,\elOv)$ is
a bi-differential of $\elOu$ and $\elOv$. A \defn{Vey Product} is a
differentiable star product where $C_r=(\tfrac{i}2{\cal P})^r/r!$ where
${\cal P}$ is the bi-differential operator given by ${\cal
P}(\elOu,\elOv) =\PB{\elOv}{\elOu}$.

\begin{countereg}
We note that, if we have a $\varepsilon$-finite OANCG, given by
$\Alg=\Algfin$, then in general, we still require an infinite
expansions in $\varepsilon$ in $\SAlg$. To see this consider the
noncommutative complex disk \cite{Klimek_Lesn1}, generated by
$\elA{z}_+,\elA{z}_-$ such that
\begin{align*}
\elA{z}_+\elA{z}_- - \elA{z}_-\elA{z}_+ &= 
\varepsilon ( 1 - \elA{z}_+\elA{z}_-) (1 - \elA{z}_-\elA{z}_+)
\end{align*}
together with the normal ordering
$\Omega(\elO{z}_-^r\elO{z}_+^s)=\elA{z}_-^r\elA{z}_+^s$.  It is easy
to see that $\Omegas^{-1}(\elA{z}_+\elA{z}_-)$ will be an infinite
expansion in $\varepsilon$. 
\end{countereg}

We have shown that an ordered algebraic noncommutative geometry gives
us a star product algebra. As mentioned in the introduction one can
define a star algebra independently simply by specifying the functions
$C_r\colon \CoM\times\CoM\mapsto\CoM$, and requiring that the star
product defined by (\ref{intro_star}) is associative.  We may now ask
whether, given such an abstractly defined star product algebra, we can
reconstruct an OANCG. For this we note the following:
\begin{jglist}
\item We construct the map $\piS\colon
\SAlg\mapsto\SOAlg=\SAlg/\Set{\varepsilon\sim0}$. Normally
$\SOAlg=\CoM$, but $\OAlg$ is an algebra of polynomials, therefore at
best we can construct an OANCG $(\Alg_1,\Omega)$ so that $\SAlg_1$,
the corresponding star algebra, is a subalgebra $\SAlg_1\subset\SAlg$.

\item An OANCG gives more information than $\SAlg$, this is given by
the immersion elements $\Set{\elSx_1,\ldots,\elSx_n}$,
$\elSx_i\in\SAlg$. These satisfy a set of immersion equations
$\QIZ_s(\underline{\elOx})=0$ where $\elOx_i=\piS(\elSx_i)$, which
defines the immersion $\man\subset\Real^n$.

\item If $\SAlg$ is constructed from $(\Alg,\Omega)$ then
$C_r(\elAx_i)\in\OAlg$ where $C_r$ is defined by
(\ref{Ord_def_Cr1}). Thus we require that $C_r(\elSx_i)$ is a
polynomial in $\elOx_k$, where $C_r\colon \SAlg\mapsto\SOAlg$ is
given by $C_r(\elSu)$ is the $\varepsilon^r$ coefficient of $\elSu$.

\item
If $\SAlg$ is constructed from $(\Alg,\Omega)$ then
$C_r(\elOx_i,\elOx_j)\in\OAlg$, thus we require that
$C_r(\elOx_i,\elOx_j)$ are polynomials in $\elOx_k$.

\item By considering the counter example above, in general it is
possible only to construct an $\varepsilon$-infinite ANCG
$\Alginf$. 

\end{jglist}
Given these conditions we can construct an OANCG. This is a precise
statement (\ref{intro_star_apprx}).

\begin{theorem}
\label{th_Ord_star_gives_A}
Given a star algebra $\SAlg$, over a manifold $\man$, and a choice of
immersion elements $\Set{\elSx_1,\ldots,\elSx_n}$, $\elSx_i\in\SAlg$
such that the set $\elOx_i=\pi(\elSx_i)$ define the immersion
$\man=\Set{\QIZ(\underline{\elOx})=0}\subset\Real^n$, and such that
$C_r(\elOx_j,\elOx_k)$ and $C_r(\elSx_i)$ are polynomials in
$(\elOx_1,\ldots,\elOx_n)$ for all $r,i,j$. Then there exists is a
unique OANCG $(\Alg_1=\Alginf_1,\man,\Omega)$ for which $\SAlg_1$,
the corresponding star algebra, is a subalgebra of $\SAlg$.
\end{theorem}

\begin{proof}
Let $\SAlg_1$ be the subalgebra of $\SAlg$ generated by star products
of $\Set{\elSx_1,\ldots,\elSx_n,\varepsilon}$.  Let $\Falg=\Falginf$ be
generated by $\Set{\elFx_1,\ldots,\elFx_n,\varepsilon}$.  Let the
function $\chi\colon \Falg\mapsto\SAlg$ be the algebraic homomorphism
satisfying
\begin{align*}
\chi(\varepsilon)=\varepsilon\,,\
\chi(\elFx_i)=\elSx_i\,,\
\chi(\elFu\elFv)=\chi(\elFu)\star\chi(\elFv)
\end{align*}
Now we show
there is a map $\Phi\colon \SAlg_1\mapsto\Falg$ such that
$\chi\circ\Phi=1_{\SAlg_1}$.

Let $\Phi'\colon \SAlg_1\mapsto\Falg$ be any map such that
$\pi\circ\q\circ\Phi'=\pi$, this can be constructed similar to proof
of theorem \ref{th_Ord_exist_Ord}. For $\elSu\in\SAlg_1$ let
$\elSu_0=\elSu$ and $\elSu_{n+1}=\elSu_n-\chi(\Phi'(\elSu_n))$.  Thus
$\elSu_n\in\SAlg$ and $\elSu_n=O(\varepsilon^n)$. Let
$\Phi(\elSu)=\sum_{n=0}^\infty \Phi'(\elSu_n)$. We can say this
converges in $\Falginf$ and it is easy to show that
$\chi\Phi(\elSu)=\elSu$.

Now let $\QC_{ij}=\elFx_i\elFx_j-\elFx_j\elFx_i- \Phi(\elSx_i\star
\elSx_j-\elSx_j\star \elSx_i)$ and let
$\QI_{s}=\Phi''(\QIZ_{s})-\Phi\chi\Phi''(\QIZ_{s})$ where
$\Phi''\colon \COFalg\mapsto\Falg$ is any map such that
$\qCO\circ\piF\circ\Phi''=1_{\COFalg}$.

The ordering is given by $\Omega=\q\circ\Phi$. This defines the
map $\Omegas\colon\SAlg_1\mapsto\Alg$.

Since $\Omegas$ is a bijective homomorphism, then $\Alg_1$ is unique.
\end{proof}

\subsection{Homomorphisms}
\label{ch_Hom}

Let $(\Alg_1,\varepsilon_1,\man_1)$ and
$(\Alg_2,\varepsilon_2,\man_2)$ be two ANCG. We say that $\Psi\colon
\Alg_1\mapsto\Alg_2$ is a \defn{homomorphism of ANCG} if $\Psi$ is
homomorphism of algebras and $\Psi(\varepsilon_1)=\varepsilon_2$.  Let
\begin{align}
&\Psi^0\colon \OAlg_1\mapsto \OAlg_2\,; 
&&
\Psi^0\circ\pi_1=\pi_2\circ\Psi
\label{Hom_Psi0}
\\
&\Psi_\star\colon \man_2\mapsto\man_1\,;
&&
\Psi^\star(\elOu) = \elOu \circ \Psi_\star\,,\qquad
\forall \elOu\in\OAlg_1
\label{Hom_def_PsiL}
\\
&\Psi^0\colon \CoMn{1}\mapsto \CoMn{2}\,;
&&
\Psi^\star(\elOu) = \elOu \circ \Psi_\star\,,\qquad
\forall \elOu\in \CoMn{1}
\end{align}
We say that $\Alg_1$ and $\Alg_2$ are \defn{isomorphic} if the map
$\Psi$ is bijective. 

\begin{theorem}
The maps $\Psi^0$, $\Psi^\star$ and $\Psi_\star$ are well
defined.  The maps $\Psi^0$ and $\Psi^\star$ are homomorphisms and
respect the Poisson structure:
\begin{align}
\Psi^\star\left(\PB{\elOu}{\elOv}\right)=
\PB{{\Psi^\star(\elOu)}}{{\Psi^\star(\elOv)}}
\end{align}
If $\Psi^0$ is surjective then $\Psi_\star$ is injective.  
If $\Psi^0$ is bijective then $\Psi_\star$ is bijective.  
Finally $\Psi$  is bijective if and only if
$\Psi^0$ is bijective 
\end{theorem}

\begin{proof}
If $\elAu\in\Alg_1$ and $\pi_1(\elAu)=0$ 
so $\elAu=\varepsilon_1 \elAu'$ for some $\elAu'\in\Alg_1$.  Thus
$\Psi(\elAu)=\Psi(\varepsilon_1\elAu')= \Psi(\varepsilon_1)\Psi(\elAu')=
\varepsilon_2\Psi(\elAu')$ so $\pi_2(\Psi(\elOu))=0$.  Thus $\Psi^0$ is
well defined, and clearly it is a homomorphism.

Let $\Set{\elAx_1,\ldots,\elAx_{n_1}}$ and
$\Set{\elAy_1,\ldots,\elAy_{n_2}}$ be the generators of $\Alg_1$ and
$\Alg_2$ respectively, and let $\elOx_i=\pi_1(\elAx_i)$ and
$\elOy_i=\pi_2(\elAy_i)$.  Given $\elOu\in\OAlg_1$ then
$\Psi^0(\elOu)\in\OAlg_2$ is a polynomial in $\elOy_i$. This includes
the functions $\Psi^0(\elOx_i)$.  Since $\Psi^0$ is a homomorphism,
then for any polynomial $f\colon \Real^{n_1}\mapsto\Cmpx$ we have
\begin{align}
f(\Psi^0(\elOx_1),\ldots,\Psi^0(\elOx_{n_1})) &=
\Psi^0(f(\elOx_1,\ldots,\elOx_{n_1}))
\label{Hom_f_Psi0}
\end{align}

Given a point $p\in\man_2$, this has coordinates
$(\elOy_1(p),\ldots,\elOy_{n_2}(p))$. The point
$\Psi_\star(p)$ has coordinates
$(\Psi^0(\elOx_1)(p),\ldots,\Psi^0(\elOx_{n_1})(p))$.
Also for each immersions equation $\QIZ_s$ defining $\Alg_1$, given in
(\ref{defANCG_QIZ}) we have  
$\QIZ_s(\Psi^0(\elOx_1),\ldots,\Psi^0(\elOx_{n_1}))=0$ 
from (\ref{Hom_f_Psi0}). So $\Psi_\star(p)\in\man_1$.  

We define $\Psi^\star$ via 
$\Psi^\star(\elOu) = \elOu \circ \Psi_\star$ for
$\elOu\in \CoMn{1}$. Since $\Psi^0(\elOx^i)$ is a polynomial in
$\elOy_i$ we can calculate all the partial derivative.   

To show that the Poisson structure is preserved, given
$\elAu,\elAv\in\Alg_1$ we have
\begin{align*}
\Psi^0(\PB{\pi_1(\elAu)}{\pi_1(\elAv})) 
&=
\Psi^0(\pi_1(\tfrac1{i\varepsilon_1}[\elAu,\elAv]))
=
\pi_2(\Psi(\tfrac1{i\varepsilon_1})[\Psi(\elAu),\Psi(\elAv)])
=
\pi_2(\tfrac1{i\varepsilon_2}[\Psi(\elAu),\Psi(\elAv)])
\\ 
&=
\PB{{\pi_2\Psi(\elAu)}}{{\pi_2\Psi(\elAv)}}
= 
\PB{{\Psi^0(\pi_1(\elAu))}}{{\Psi^0(\pi_1(\elAv))}}
\end{align*}
And since the Poisson bracket is defined by a bi-vector then 
$\Psi^\star$ must also respect the Poisson structure.

Let $\Psi^0$ be surjective and $p,q\in\man_2$ such that $p\ne q$. Then
there exists a function $\elOu\in \OAlg_2$ such that
$\elOu(p)\ne\elOu(q)$. Since $\Psi^0$ is surjective there exists a
function $\elOv\in \OAlg_1$ such that $\Psi^0(\elOv)=\elOu$. Thus
$\Psi^0(\elOv)(p)\ne\Psi^0(\elOv)(q)$, giving
$\elOv(\Psi_\star(p))\ne\elOv(\Psi_\star(q))$, and hence
$\Psi_\star(p)\ne\Psi_\star(q)$. Thus $\Psi_\star$ is injective.

If $\Psi^0$ is bijective then we have
$\Psi^{0-1}\colon \OAlg_2\mapsto \OAlg_1$ and hence 
$\Psi_\star^{-1}\colon \man_1\mapsto\man_2$. Given $\elOu\in
\OAlg_1$, $p\in\man_1$ we have
\begin{align*}
\elOu(\Psi_\star^{-1}(\Psi_\star(p))) =
\Psi^{0-1}(\elOu)(\Psi_\star(p)) =
\Psi^{0-1}(\Psi^\star(\elOu))(p)=\elOu(p)
\end{align*}
Since this is for all $\elOu$, then $\Psi_\star$ is bijective.

If $\Psi$ is bijective then we can define $\Psi^{0-1}$ via
$\Psi^{0-1}\circ\pi_2=\pi_1\circ\Psi^{-1}$, and this satisfies
$\Psi^{0-1}\Psi^{0}=1_{\OAlg_1}$ and $\Psi^{0}\Psi^{0-1}=1_{\OAlg_2}$.

If $\Psi^0$ is bijective then we first show that $\Psi$ is
injective. Let $\elAu\in\Alg_1$ such that $\Psi(\elAu)=0$. Then
$\pi_2\circ\Psi(\elAu)=\Psi^0\circ\pi_1(\elAu)=0$. Since $\Psi^0$ is
injective $\pi_1(\elAu)=0$. Thus $\elAu=\varepsilon_1 \elAu_1$. So
$0=\Psi(\elAu)=\varepsilon_2\Psi(\elAu_1)$. So
$\Psi(\elAu_1)=0$. Repeating this process shows $\elAu=0$.

If $\Psi^0$ is bijective then we show, by construction, that $\Psi$ is
surjective. Choose any ordering $\Omega_2\colon
\OAlg_2\mapsto\Alg_2$. Let $\elAv_n\in\Alg_2$, $n=0,1,\ldots$ be
defined inductively via
\begin{align*}
\elAv_0 &= \Omega_2\circ\Psi^{0-1}\circ\pi_1(\elAu) 
\\
\elAv_{n+1} &=  \elAv_n + 
\varepsilon_2^n\Omega_2\circ\Psi^{0-1}\circ\pi_1
\left(\frac{\elAu-\Psi(\elAv_n)}{\varepsilon_1^n}\right)
\end{align*}
Clearly $\elAv_{n+1}-\elAv_{n}=O(\varepsilon_1^n)$ so $\elAv_n$
converge to $\elAv_n\to \elAv\in\Alg_1$. Also
$\elAu-\Psi(\elAv_n)=O(\varepsilon_1^n)$ so $\elAu=\Psi(\elAv)$.
\end{proof}

\vspace{1 em}

If $(\Alg_1,\man_1,\Omega_1)$ and $(\Alg_2,\man_2,\Omega_2)$ are two
OANCG then we say $\Psi\colon \Alg_2\mapsto\Alg_1$ is an \defn{ANCG
homeomorphism which respects the ordering} if
\begin{align}
\Psi\circ\Omega_2=\Omega_1\circ\Psi^0
\label{HomO_def}
\end{align}
This gives the following theorem.
\begin{theorem}
\label{th_Hom_Ord_star}
If $\Psi\colon (\Alg_1,\Omega_1)\mapsto(\Alg_2,\Omega_2)$ is a ANCG
homeomorphism which respects ordering then
\begin{align}
\Psi\circ\Omegas_2=\Omegas_1\circ\Psi^0
\label{Hom_PsiOm}
\end{align}
and
\begin{align}
\Psi^0(\elOu\star_1 \elOv)= \Psi^0(\elOu)\star_2 \Psi^0(\elOv)
\end{align}
for $\elOu,\elOv\in\OAlg_1$, where $\star_1$ and $\star_2$ are the
star products corresponding to $\SAlg_1$ and $\SAlg_2$ respectively,
and $\Psi^0:\SAlg_1\mapsto\SAlg_2$ is defined via linear extension,
with $\Psi^0(\varepsilon_1)=\Psi^0(\varepsilon_2)$.
\end{theorem}
\begin{proof}
Since (\ref{HomO_def}) is linear, then we can extend $\Psi$ and
$\Omega$ to $\SAlg$, thus (\ref{Hom_PsiOm}).
Now
\begin{align*}
\Psi^0(\elOu\star_1 \elOv) 
&= 
\Psi^0 \circ \Omegas^{-1}_1 (\Omega_1(\elOu)\Omega_1(\elOv)) 
=
\Omegas^{-1}_2 \circ \Psi (\Omega_1(\elOu)\Omega_1(\elOv)) 
=
\Omegas^{-1}_2  (\Psi\circ\Omega_1(\elOu) \Psi\circ\Omega_1(\elOv)) 
\\
&=
\Omegas^{-1}_2  (\Omega_2\circ\Psi^0(\elOu) \Omega_2\circ\Psi^0(\elOv)) 
=
\Psi^0(\elOu)\star_2\Psi^0(\elOv)
\end{align*}
\end{proof}


\subsection{Representations and Trace}
\label{ch_Rep}

An additional structure that an ANCG may have is a representation or
matrix representation.  This is independent of whether or not an
ordering is specified.

A \defn{representations} of $\Alg$ over the Hilbert space $\Hil$ is a
homomorphism
\begin{align}
\varphi\colon \Alg\hookrightarrow L(\Hil)\,;
\qquad
\varphi(\varepsilon)=\varepsilon_\infty \in\Real
\end{align}
Here $L(\Hil)$ is the space of linear (but not necessarily bounded)
operators on $\Hil$.  This representation is unitary if
$\varphi(\elAu^\dagger)=\varphi(\elAu)^\dagger$ where
$\varphi(\elAu)^\dagger$ is the adjoint with respect to the inner
product on $\Hil$. 

Clearly if $\man$ is compact and $\varepsilon_\infty =0$ there is a
natural unitary representation with $\Hil=L^2(\man)$ as
$\varphi(\elAu)f=\pi(\elAu)f$ with $f\in L^2(\man)$.  If
$\varepsilon_\infty \ne0$ then a prerequisite for the existence of a
representation is that $\Alg=\Algfin$. This is because the element
$\sum_{r=0}^\infty \varepsilon^r r!\in\Alginf$, and this does not have
an image under $\varphi$.

We say there is a \defn{matrix approximation} of $\Alg$ if there exists
a sequence of $\varepsilon_N\in\Real$ with $\varepsilon_N\ne0$, and
$\varepsilon_N\to0$ as $N\to\infty$, such that
\begin{align}
\varphi_N\colon \Alg\mapsto L(\Cmpx^N)=M_N(\Cmpx)\,;
\qquad \varphi_N(\varepsilon)=\varepsilon_N
\end{align}
Given an ANCG, it is not a trivial matter deciding whether there
exists a unitary representation. In section \ref{ch_eg} we give
a number of examples of ANCGs with representation. Here we give an ANCG
for a compact manifold,
which does not possess a unitary matrix representation.
\begin{countereg}
This counter example is given by tensoring two copies of the
noncommutative torus given in section \ref{ch_egT2}.
Let $\man=T^4$, and $\Alg$ be generated by
$\Set{\varepsilon,\elAu_1,\elAu_2,\elAv_1,\elAv_2,
\elAu_1^{-1},\elAu_2^{-1},\elAv_1^{-1},\elAv_2^{-1}}$.
These obey
\begin{align}
&\elAu_1^r\elAv_1^s=e^{i\varepsilon rs}\elAv_1^s\elAv_1^r\,,
&
&\elAu_2^r\elAv_2^s=e^{i\varepsilon\alpha rs}\elAv_2^s\elAv_2^r\,,
&
&\elAu_i^\dagger=\elAu_i^{-1}\,,
&
&\elAv_i^\dagger=\elAv_i^{-1}
\end{align}
where $\alpha\in\Real\backslash\Rationals$ is an irrational number,
$i=1,2$ and all other commutators are zero.
\begin{lemma}
The above ANCG does not have a unitary matrix representation. 
\end{lemma}
\begin{proof}
Let us assume we have a representation $\varphi_N$.  Let
$\lambda_1,\ldots,\lambda_N$ be the eigenvalues of
$\varphi_N(\elAu_1)$. Then $\sum_i \lambda_i^r=
\Tr(\varphi_N(\elAu_1^r))$. This implies that there must exist an
$r_1>0$ such that $\Tr(\varphi_N(\elAu_1^{r_1}))\ne0$.  By looking at
the trace of $\elAv_1\elAu^r_1\elAv^{-1}_1$ we can show that
$(1-e^{ir\varepsilon_N})\Tr(\varphi_N(\elAu_1^r))=0$.  Hence
$e^{ir_1\varepsilon_N}=1$. Likewise $e^{ir_2\varepsilon_N\alpha}=1$,
for another integer $r_2>0$. This is impossible since $\alpha$ is not
rational.
\end{proof}
\end{countereg}

Given two ANCG, $\Alg_1$ and $\Alg_2$, with a homomorphism $\Psi\colon
\Alg_1\mapsto\Alg_2$. If both ANCG have matrix representations
respectively given by $\varphi_N^{(1)}$ and $\varphi_N^{(2)}$, then
this induces a matrix homomorphism $\Psi_N$ for each $N$ given by
\begin{align}
\Psi_N\colon M_N(\Cmpx)\mapsto M_N(\Cmpx)\,;\qquad
\Psi_N\circ\varphi_N^{(1)} =
\varphi_N^{(2)}\circ\Psi
\end{align}
Alternatively if only $\Alg_2$ has a representation $\varphi_N^{(2)}$
then we can induce a representation of $\Alg_1$
\begin{align}
\varphi_N^{(1)} = \varphi_N^{(2)}\circ\Psi
\end{align}
However the representation generated in the way will not be
surjective, unless $\Psi$ is an isomorphism (see below). 
In this case we have the following trivial lemma.
\begin{lemma}
If $\Alg_1$ and $\Alg_2$ are isomorphic ANCG and $\varphi_N^{(2)}$ is a
surjective matrix representation, then $\varphi_N^{(1)}$ is also a
surjective matrix representation.
\end{lemma}

\vspace{1 em}

If ($\Alg,\man$) has a matrix approximation we define the the
\defn{trace function} as the map
\begin{align}
\Trn{N}\colon \Alg\mapsto\Cmpx\,; 
\qquad
\Trn{N}(\elAu)=\tfrac1{N}\Tr(\varphi_N(\elAu))
\label{Rep_TrN}
\end{align}
where $\Tr\colon M_N(\Cmpx)\mapsto\Real$ is the matrix trace.
In general the matrix trace is dependent on the choice of matrix
approximations however we do have the following trivial lemma
\begin{lemma}
If $U_N\in GL_N(\Cmpx)$ and $\varphi_N$ is a matrix approximation of
$(\Alg,\man)$ then $\varphi'_N=U_N\varphi_N U_N^{-1}$ defines
another matrix approximation of $\Alg$. In this case
$\Trn{N}(\elAu)=\Trn{N}'(\elAu)$.
\end{lemma}

\begin{theorem}
\label{thm_int}
If $(\Alg,\man)$ is symplectic, $\man$ is compact  
and $\Trn{N}$ exists for all $N$ and 
$\lim_{N\to\infty}\Trn{N}(\elAu)$ converges for all $\elAu\in\Alg$
then
\begin{align}
\lim_{N\to\infty}\Trn{N}(\elAu)
=
\frac1{|\man|}\int_\man \pi(\elAu)\omega^r
\end{align}
where $|\man|=\int_\man\omega$, and $\dim(\man)=2r$.
\end{theorem}

\begin{proof}
Choose some ordering on $\Alg$. Define $\TrO\colon \OAlg\mapsto\Cmpx$
as $\TrO(\elOu)=\lim_{N\to\infty}\Trn{N}(\Omega(\elOu))$. Since $\TrO$
is linear on $\OAlg$ we can write
\begin{align}
\TrO(\elOu)=\int_\man \elOu W \omega^r
\end{align}
for some (distributional) weight function $W$ on $\man$.  

Let $(p_1,\ldots,p_r,q_1,\ldots,q_r)$ be conjugate coordinates on a
patch of $\man$, and let $x_{2s}=p_s$ and $x_{2s+1}=q_s$. Let
$\elAu_s\in\Alg$ and $\elOu_s=\pi(\elAu_s)$, for $s=1,\ldots,2r$
\begin{align*}
\lefteqn{
\sum_{\sigma\in S_{2r}} \epsilon(\sigma)
\PB{\elOu_{\sigma(1)}}{\elOu_{\sigma(2)}}
\cdots
\PB{\elOu_{\sigma(2r-1)}}{\elOu_{\sigma(2r)}}
\omega^r}
\qquad\qquad
&
\\
&=
r!
\sum_{\sigma\in S_{2r}} \epsilon(\sigma)
\sum_{i_1,\ldots i_r=1}^{r}
\left(
\frac{\partial \elOu_{\sigma(1)}}{\partial p_{i_1}}
\frac{\partial \elOu_{\sigma(2)}}{\partial q_{i_1}}
-
\frac{\partial \elOu_{\sigma(2)}}{\partial p_{i_1}}
\frac{\partial \elOu_{\sigma(1)}}{\partial q_{i_1}}
\right)
\cdots
\times
\\
&\qquad\qquad\qquad\qquad
\left(
\frac{\partial \elOu_{\sigma(2r-1)}}{\partial p_{i_r}}
\frac{\partial \elOu_{\sigma(2r)}}{\partial q_{i_r}}
-
\frac{\partial \elOu_{\sigma(2r)}}{\partial p_{i_r}}
\frac{\partial \elOu_{\sigma(2r-1)}}{\partial q_{i_r}}
\right)
dp_1\wedge dq_1\wedge \cdots \wedge dp_r\wedge dq_r
\\
\qquad
&=
2^r r!
\sum_{\sigma\in S_{2r}} \epsilon(\sigma)
\sum_{i_1,\ldots i_r=1}^{r}
\frac{\partial \elOu_{\sigma(1)}}{\partial p_{i_1}}
\frac{\partial \elOu_{\sigma(2)}}{\partial q_{i_1}}
\cdots
\frac{\partial \elOu_{\sigma(2r-1)}}{\partial p_{i_r}}
\frac{\partial \elOu_{\sigma(2r)}}{\partial q_{i_r}}
dp_1\wedge dq_1\wedge \cdots \wedge dp_r\wedge dq_r
\\
&=
2^r r!
\sum_{\sigma\in S_{2r}} \epsilon(\sigma)
\sum_{\tau\in S_r}
\frac{\partial \elOu_{\sigma(1)}}{\partial p_{\tau(1)}}
\frac{\partial \elOu_{\sigma(2)}}{\partial q_{\tau(1)}}
\cdots
\frac{\partial \elOu_{\sigma(2r-1)}}{\partial p_{\tau(r)}}
\frac{\partial \elOu_{\sigma(2r)}}{\partial q_{\tau(r)}}
dp_1\wedge dq_1\wedge \cdots \wedge dp_r\wedge dq_r
\\
&=
2^r (r!)^2
\sum_{\sigma\in S_{2r}} \epsilon(\sigma)
\frac{\partial \elOu_{\sigma(1)}}{\partial p_1}
\frac{\partial \elOu_{\sigma(2)}}{\partial q_1}
\cdots
\frac{\partial \elOu_{\sigma(2r-1)}}{\partial p_r}
\frac{\partial \elOu_{\sigma(2r)}}{\partial q_r}
dp_1\wedge dq_1\wedge \cdots \wedge dp_r\wedge dq_r
\\
&=
2^r (r!)^2
\det_{ij}\left(
\frac{\partial \elOu_i}{\partial x_j}
\right)
dx_1\wedge dx_2\wedge \cdots \wedge dx_{2r-1}\wedge dx_{2r}
\\
&=
2^r (r!)^2
d\elOu_1\wedge d\elOu_2\wedge \cdots \wedge d\elOu_{2r}
\end{align*}
where $S_{r}$ is the set of permutations, and $\epsilon(\sigma)$ is
the signature of the permutation. 
However
\begin{align*}
\lefteqn{
\int_\man \omega^r
\sum_{\sigma\in S_{2r}} \epsilon(\sigma)
\PB{\elOu_{\sigma(1)}}{\elOu_{\sigma(2)}}
\cdots
\PB{\elOu_{\sigma(2r-1)}}{\elOu_{\sigma(2r)}}
}\qquad\qquad
&
\\
&=
\lim_{N\to\infty}
\sum_{\sigma\in S_{2r}} \epsilon(\sigma)
\Trn{N}
\left(
(i\varepsilon)^{-r}
[\elAu_{\sigma(1)},\elAu_{\sigma(2)}]
\cdots
[\elAu_{\sigma(2r-1)},\elAu_{\sigma(2r)}]
\right)
\\
&=
\lim_{N\to\infty}
2^r (i\varepsilon_N)^{-r}
\sum_{\sigma\in S_{2r}} \epsilon(\sigma)
\Trn{N}
\left(
\elAu_{\sigma(1)}\elAu_{\sigma(2)}
\cdots\elAu_{\sigma(2r)}
\right)
\\
&=
\lim_{N\to\infty}
2^r (i\varepsilon_N)^{-r}
\sum_{\sigma\in S_{2r}} \epsilon(\sigma)
\Trn{N}
\left(
\elAu_{\sigma(2)}
\cdots\elAu_{\sigma(2r)}\elAu_{\sigma(1)}
\right)
=0
\end{align*}
since it is an odd permutation. Thus we have
\begin{align}
\int_\man W d\elOu_1\wedge d\elOu_2\wedge \cdots \wedge d\elOu_{2r}=0
\end{align}
for all (algebraic) functions $\elOu_s$ on $\man$. Integration by
parts gives
\begin{align*}
\int_\man \elOu_{1} 
dW\wedge d\elOu_2\wedge \cdots \wedge d\elOu_{2r}=0
\end{align*}
By considering a sequence of $\elOu_1$ we can let $\elOu_1$ be the
characteristic function on some subset $U_1\in\man$. This implies
\begin{align*}
0=
\int_{U_1}  
d W\wedge d\elOu_2\wedge \cdots \wedge d\elOu_{2r}
=
-
\int_{\partial U_1}  
\elOu_2 d W\wedge d\elOu_3\wedge \cdots \wedge d\elOu_{2r}
\end{align*}
Repeat this process until we are left with $\int_a^b dW=0$, for two
points $a,b\in\man$. This implies $W(a)=W(b)$. Thus $W$ is a constant,
whose value is given by $\Trn{N}(1)=1$.
\end{proof}

If $\Trn{N}$ exists, we can define an inner product form:
\begin{align}
\innerprod{\bullet,\bullet}_N \colon  
\Alg\times\Alg
\mapsto \Cmpx 
\qquad
\innerprod{\elAu,\elAv}_N = \Trn{N}(\elAu^\dagger\elAv) 
\end{align}
This obeys
$\innerprod{[\elAu,\elAv],\elAw}_N=\innerprod{\elAu,[\elAv,\elAw]}_N$.

We say $\Trn{N}$ is \defn{analytic} if we can define a function
$\TrA\colon \Alg\mapsto C^\omega(\Real)$, such that
$\TrA(\elAu)(\varepsilon_N)=\Trn{N}(\elAu)$ for all $N$.  For example
the trace function on the noncommutative sphere and surface of
rotation is analytic whilst the trace function on the noncommutative
torus is not analytic, (see sections
\ref{ch_egSR},\ref{ch_egS2},\ref{ch_egSR}). If $\TrA$ exists then we
can define the sesquilinear form
\begin{align}
\innerprod{\bullet,\bullet}\colon\Alg\times\Alg\mapsto
C^\omega(\Real)\,;
\qquad
\innerprod{\elAu,\elAv}(\varepsilon)=
\TrA(\elAu^\dagger\elAu)(\varepsilon) 
\end{align}
which satisfies
$\innerprod{\elAu,\elAv}(\varepsilon_N)=\innerprod{\elAu,\elAv}_N$.
Although $\innerprod{\bullet,\bullet}_N$ is an inner product,
$\innerprod{\bullet,\bullet}$ is, in general,  not positive definite
for all $\varepsilon\in\Real$.

If $(\Alg,\Omega)$ is an OANCG  then we say $\Trn{N}$ is
compatible with $\Omega$ if $\Trn{N}(\Omega(\elOu))$ is
independent of $N$. One example is the sphere with the
Wick-like ordering given by (\ref{egS2_OmS2}).


\subsection{The Heisenberg Algebra and Coordinate Charts}
\label{ch_Heis}

For $r,s\in\Intg$, $r\ge1$, $s\ge0$, we call the \defn{Heisenberg
algebra} $\Heis_{2r,s}$ the algebra generated by
$\Set{\elAp_1,\ldots,\elAp_r,\elAq_1,\ldots,\elAq_r,
\elAy_1,\ldots,\elAy_s,\varepsilon}$ with the only nonzero $\QC_{ij}$
given by $[\elAp_i,\elAq_j]=i\varepsilon\delta_{ij}$ and with no
immersions equation $\QI_{t}$ so $\dimman=2r+s$. Clearly this is a
ANCG for the manifold $\Real^{2r+s}$. Thus $\Heis_{2r,s}^0$ is the
algebra of polynomials on $\Real^{2r+s}$. Each $\elAy_i$ is in the
centre of $\Heis_{2r,s}$ so the corresponding symplectic leaves of
$\Real^{2r+s}$ are given by $\elOy_i=\text{constant}$. If $s=0$ then
we define $\Heis_{2r}=\Heis_{2r,s}$.

Two orderings on $\Heis_{2r,s}$ are commonly considered,  the
\defn{Wick} ordering and the \defn{normal} ordering.
The Wick ordering is unique and is given by
\begin{equation}
\begin{aligned}
\Omega_W(
\elOp_1^{i_1}\cdots\elOp_r^{i_r}
\elOq_1^{j_1}\cdots\elOq_r^{j_r}
\elOy_1^{k_1}\cdots\elOy_s^{k_s})
=&
\text{ the correctly normalised sum of all symmetric} 
\\
&\text{permutations of }
\elAp_1^{i_1}\cdots\elAp_r^{i_r}
\elAq_1^{j_1}\cdots\elAq_r^{j_r}
\elAy_1^{k_1}\cdots\elAy_s^{k_s}
\end{aligned}
\label{Heis_Om_W}
\end{equation}
where correctly normalised means that
$\pi\circ\Omega_W=1_{\Heis_{2r,s}^0}$. The number of terms in
the symmetric sum is given by
\begin{align}
\frac{(i_1+\cdots+i_r+j_1+\cdots+j_r+k_1+\cdots+k_s)!}
{i_1!\cdots i_r! j_1! \cdots j_r! k_1!\cdots k_s!}
\end{align}
so we must divide by this quantity. 

The normal ordering depends on the choices of an ordering on the
generators of $\Heis_{2r,s}$. It is related to the time ordering in
quantum field theory. The choice we will use here is to place 
$\elAp_i$ before $\elAq_i$ thus
\begin{align}
&\Omega_N(
\elOp_1^{i_1}\cdots\elOp_r^{i_r}
\elOq_1^{j_1}\cdots\elOq_r^{j_r}
\elOy_1^{k_1}\cdots\elOy_s^{k_s})
=
\elAp_1^{i_1}\cdots\elAp_r^{i_r}
\elAq_1^{j_1}\cdots\elAq_r^{j_r}
\elAy_1^{k_1}\cdots\elAy_s^{k_s}
\label{Heis_Om_N}
\end{align}

For the Heisenberg plane $\Heis_2$ let the Wick basis elements
$\elA{S}(a,b)=\Omega_W(\elOp^a\elOq^b)$ and the normal basis elements be
$\elA{N}(a,b)=\elAp^a\elAq^b=\Omega_N(\elOp^a\elOq^b)$. 

\begin{theorem}
The Wick basis elements and normal basis elements are related by 
\label{th_Heis_SN}
\begin{align}
&
\elA{S}(a,b)=
\sum_{r=0}^{\min(a,b)} 
\frac{(-\tfrac12i\varepsilon)^r}{r!}
\frac{a!}{(a-r)!}\frac{b!}{(b-r)!}
\elA{N}(a,b)
\label{Heis_S2N}
\\
&
\elA{N}(a,b)=
\sum_{r=0}^{\min(a,b)} 
\frac{(\tfrac12i\varepsilon)^r}{r!}
\frac{a!}{(a-r)!}\frac{b!}{(b-r)!}
\elA{S}(a,b)
\label{Heis_N2S}
\end{align}
The product of two basis elements are given by
\begin{align}
&
\elA{N}(a,b)\elA{N}(c,d)=
\sum_{r=0}^{\min(b,c)} 
\frac{(-i\varepsilon)^r}{r!}
\frac{b!}{(b-r)!}
\frac{c!}{(c-r)!}
\elA{N}(a+c-r,b+d-r)
\label{Heis_NxN}
\\
&
\elA{S}(a,b)\elA{S}(c,d)=
\sum_{r=0}
\frac{(i\varepsilon)^r}{r!}
\elA{S}(a+c-r,b+d-r)
\sum_{s=0}^n 
\frac{(-1)^{r-s}}{s!(r-s)!}
\frac{a!}{(a-r+s)!}
\frac{b!}{(b-s)!}
\frac{c!}{(c-r+s)!}
\frac{d!}{(d-s)!}
\label{Heis_SxS}
\end{align}
\end{theorem}

\begin{proof}
First note
\begin{align*}
\elAp\elA{S}(a,b)+\elA{S}(a,b)\elAp
&=
2\elA{S}(a+1,b)\,,
&
\elAq\elA{S}(a,b)+\elA{S}(a,b)\elAq
&=
2\elA{S}(a,b+1)\,,
\\
{}[\elAp,\elA{S}(a,b)]
&=
i\varepsilon b \elA{S}(a,b-1)\,,
&
{}[\elAq,\elA{S}(a,b)]
&=
-i\varepsilon a\elA{S}(a-1,b)\,.
\end{align*}
These are given in \cite[appendex]{JG-TS2}. Also
\begin{align*}
\elAp\elA{N}(a,b)+\elA{N}(a,b)\elAp
&=
2\elA{N}(a+1,b)+ \tfrac12 b i \varepsilon 
\elA{N}(a,b-1)
\end{align*}
Thus (\ref{Heis_N2S}) follows from induction on $a$, and
(\ref{Heis_S2N}) is its inverse.

Equation (\ref{Heis_NxN}) follows from induction on $b$.
For (\ref{Heis_SxS}) expand $\elAp\elA{S}(a,b)+\elA{S}(a,b)\elAp$ and
$\elAq\elA{S}(a,b)+\elA{S}(a,b)\elAq$. Then (\ref{Heis_SxS}) follows
from induction on $a$ and $b$.
\end{proof}

\begin{theorem}
The star product on $\Heis_{2r,s}$ with the Wick ordering is the Vey
product.
\begin{align}
\elOu\star_W\elOv=\exp(\tfrac12i\varepsilon{\cal P})(\elOu,\elOv)
\label{Heis_starW}
\end{align}
where ${\cal P}$ is the Poisson operator given by ${\cal
P}(\elOu,\elOv)=\PB{\elOu}{\elOv}$. That is
\begin{align}
{\cal P}=
\sum_i\left(
\frac{\partial_1}{\partial \elOp_i}
\frac{\partial_2}{\partial \elOq_i}
-
\frac{\partial_2}{\partial \elOp_i}
\frac{\partial_1}{\partial \elOq_i}
\right)
\end{align}
where the subscript $1,2$ refer to differentiation with respect to the
first and second variable.

The star product on $\Heis_{2r,s}$ with the Normal ordering is
\begin{align}
\elOu\star_N\elOv=\exp(-i\varepsilon{\cal N})(\elOu,\elOv)
\label{Heis_starN}
\qquad\text{where}\qquad
{\cal N}=
\sum_i
\frac{\partial_2}{\partial \elOp_i}
\frac{\partial_1}{\partial \elOq_i}
\end{align}
\end{theorem}
\begin{proof}
To show this is true for $\Heis_2$ simply substitute 
$\elA{N}(a,b)$ into (\ref{Heis_starN}) and $\elA{S}(a,b)$ into
(\ref{Heis_starN}) to obtain the corresponding product formulae.
The results naturally extend for $\Heis_{2r,s}$.
\end{proof}

\vspace{1 em}

In order to interpret $\Heis_{2r,s}$ as a coordinate basis we need to
enlarge it to include certain analytic functions of the generators.

Let $\underline{a},\underline{b}\in
\Big(\Real\union\Set{\pm\infty}\Big)^{2r+s}$ such that $a_i<b_i$.  Let
$\Heis_{2r,s}(\underline{a},\underline{b})$ be the algebra generated
by $\Set{f_i(\elAp_i),g_i(\elAq_i),h_i(\elAy_i),\varepsilon}$ (with
infinite sums of $\varepsilon$) where $f_i\in C^\omega(a_i,b_i)$,
$g_i\in C^\omega(a_i+r,b_i+r)$ and $h_i\in C^\omega(a_i+2r,b_i+2r)$,
and where $C^\omega(a_i,b_i)$ is the space of analytic functions on
$\Set{x|a_i< x< b_i}$.  The following lemma shows that
$\Heis_{2r,s}(\underline{a},\underline{b})$ is an algebra.

\begin{lemma}
Every element of $\Heis_{2r,s}(\underline{a},\underline{b})$ 
may be written in the form 
\begin{align}
\elAu=\sum_{t=0}^\infty \varepsilon^t \elAu_t
\end{align}
where $\elAu_t$ is a finite sum of terms of the form
\begin{align}
f_1(\elAp_1)\cdots f_r(\elAp_r)
g_1(\elAq_1)\cdots g_r(\elAq_r)
h_1(\elAy_1)\cdots h_s(\elAy_s)
\label{Heis_gen_el}
\end{align}
\end{lemma}
\begin{proof}
The formula for the normal star product $\Omega_N$ extends
naturally to the elements of
$\Heis^0_{2r,s}(\underline{a},\underline{b})$. Thus
\begin{align}
g_i(\elAq_i)f_i(\elAp_i)
=
\sum_{r=0}^\infty 
\frac{(-i\varepsilon)^r}{r!}f_i^{(r)}(\elAp_i) g_i^{(r)}(\elAq_i)
\end{align}
Hence result.
\end{proof}

Given $\Alg$ with $\dim(\man)=\dimman={2r+s}$ we say
there exists a \defn{Heisenberg coordinate chart} of $\Alg$ if
there exists an injective homeomorphism of ANCG
$\Psi\colon \Alg\mapsto\Heis_{2r,s}(\underline{a},\underline{b})$.

\begin{lemma}
If $\Alg$ is symplectic and $\Heis_{2r,s}$ is a coordinate chart for
$\CAlg$ then $\Heis_{2r}$ is a coordinate chart for $\Alg$. And the
local immersions relations are 
\begin{align}
\elAy_i=0
\end{align}
\end{lemma}
\begin{proof}
Trivial.
\end{proof}


\subsection{Quantum Groups}
\label{ch_QG}

We can give many ANCGs a quantum group structure as a result of the two
following theorems. 
\begin{theorem}
Let $\Psi:\Alg_1\mapsto\Alg_2$ be a isomorphism of ANCG, and let
$\Alg_1$ be a quantum group with coproduct $\Delta_1$, counit
$\epsilon_1$ and antipode $S_1$, then $\Alg_2$ is also a quantum group
with
\begin{align}
&\Delta_2=(\Psi\otimes\Psi)\circ\Delta_1\circ\Psi^{-1}\,,
&
&\epsilon_2=\epsilon_1\circ\Psi^{-1}\,,
&
& S_2 = \Psi\circ S_1 \circ \Psi^{-1}\,.
\label{QG_Hom}
\end{align}
\end{theorem}
\begin{proof}
Simply go though all the axioms of a quantum group.
\end{proof}

\begin{theorem}
The Heisenberg ANCG is a Quantum Group. 
\begin{equation}
\begin{matrix}
\Delta(1)=1\otimes 1 \,,
&
\qquad\quad
\Delta(\varepsilon)=\varepsilon \otimes 1 + 1 \otimes \varepsilon\,,
\qquad\quad
&
\Delta(\elAx)=\elAx\otimes 1 + 1\otimes \elAx\,,
\\
\epsilon(1)=1  \,,
&
\epsilon(\varepsilon)=0   \,,
&
\epsilon(\elAx)=0   \,,
\\
 S(1) = 1 \,,
& 
 S(\varepsilon) = -\varepsilon \,,
& 
 S(\elAx) = -\elAx \,,
\\
\multicolumn{3}{c}{
\forall \elAx\in\Set{\elAp_1,\ldots,\elAp_r,\elAq_1,\ldots,\elAq_r,
\elAy_1,\ldots,\elAy_s}}
\end{matrix}
\label{QG_Heis}
\end{equation}
\end{theorem}
\begin{proof}
Simply go though all the axioms of a quantum group.
\end{proof}

We can use these theorems to give a quantum group structure to ANCG
with coordinate charts. This will be used in the examples of the
noncommutative torus and surface of rotation.

\subsection{Generating a New ANCG by Use of a Homomorphism}
\label{ch_new}

Assume we have an ANCG $(\Alg_1,\man_1,\Omega_1)$, where
$\man_1\subset\Real^{n_1}$, a second manifold
$\man_2\subset\Real^{n_2}$, and an analytic bijective diffeomorphism
$\Psi_\star:\man_2\mapsto\man_1$.  We can ask whether we can generate
an OANCG $(\Alg_2,\man_2,\Omega_2)$ and a isomorphisms $\Psi\colon
\Alg_1\mapsto\Alg_2$ which respects ordering. This is important for
the application later on when we wish to construct an OANCG on a
general manifold. Unfortunately, in general, this is not possible.
However, if $\star_1$ is differentiable, we use this to define the
algebra $\SAlg_2$ via
\begin{align}
\elSu\star_2 \elSv= 
\Psi^{\star}(\Psi^{\star-1}(\elSu)\star_1 \Psi^{\star-1}(\elSv))
\end{align}
We also define the immersions elements
$\Set{\elSy_1,\ldots,\elSy_{n_2}}$, $\elSy_i\in\SAlg_2$ as the
coordinates of $\Real^{n_2}$.  However we cannot use theorem
\ref{th_Ord_star_gives_A}, because we can not guarantee that
$C_r^{(2)}(\elSy_i,\elSy_j)$ is a polynomial.

Alternatively, if $\Alg_1$ has a matrix representation, we can use
that.  Let us assume that $(\Alg_2,\Omega_2)$ does exist, and let
$\elAx_i\in\Alg_1$ and $\elAy_i\in\Alg_2$ be the corresponding bases.
Then clearly $\varphi_N^{(2)}\circ\Omega_2=
\varphi_N^{(1)}\circ\Omega_1\circ\Psi^{\star-1}$. So we have the
matrix
$Y^{(N)}_i=\varphi_N^{(2)}(\elAy_i)=
\varphi_N^{(2)}\circ\Omega_2(\elOy_i)$.
Thus
\begin{align}
Y^{(N)}_i = \varphi_N^{(1)}\circ\Omega_1\circ\Psi^{\star-1}(y_i)
\label{iso_Yi}
\end{align}
However we can define $Y^{(N)}_i\in M_N(\Cmpx)$ using (\ref{iso_Yi})
even if $\Alg_2$ does not exist.


\subsection{Geometric Properties of Surfaces}
\label{ch_Geo}

For many applications, especially gravity, we are interested in the
geometric structure of $\man$, arising from a metric. Of course we are
completely free to choose any metric on $\man$. However since we have
the embedding $\colon \man\hookrightarrow\Real^n$, we shall choose the
metric $\man$ to be the pullback of the Euclidean metric on
$\Real^n$. Let $\sharp:T^\star\man\mapsto T\man$ be the metric dual
given by $\xi(X)=g(\xi^\sharp,X)$. For this chapter we shall only
consider two dimensional surfaces immersed in $\Real^n$.

\begin{theorem}
Let $\man$ be a surface embedded in $\Real^n$ and let $(p,q)$ be
conjugate coordinates with $\PB{p}{q}=1$ on a patch $U\subset\man$.
The metric can be given solely in terms of the Poisson structure and
the functions $x_i,p,q\colon U\mapsto\Real$
\begin{align}
g({du}^\sharp,{dv}^\sharp)=
\frac1C  \sum_i \PB{x_i}{u}\PB{x_i}{v}
\end{align}
where $C\colon U\mapsto\Real$ is given by
\begin{align}
C=\sum_{ij} \PB{p}{x_i}\PB{q}{x_j}\PB{x_j}{x_i}
\label{Geo_C}
\end{align}
\end{theorem}

\begin{proof}
This is basic manipulation
\begin{align*}
g &= \sum_i dx_i\otimes dx_i 
= \sum_i \left(
\left(\frac{\partial x_i}{\partial p}\right)^2 dp\otimes dp +
\left(\frac{\partial x_i}{\partial q}\right)^2 dq\otimes dq +
\frac{\partial x_i}{\partial p}\frac{\partial x_i}{\partial q} 
(dp\otimes dq +dq\otimes dp )
\right)
\end{align*}
Inverting this gives
\begin{align*}
{g}(du^\sharp,dv^\sharp) &=\frac1C \sum_i \left(
\left(\frac{\partial x_i}{\partial q}\right)^2 
\frac{\partial u}{\partial p}\frac{\partial v}{\partial p} 
+
\left(\frac{\partial x_i}{\partial p}\right)^2 
\frac{\partial u}{\partial q}\frac{\partial v}{\partial q} 
-
\frac{\partial x_i}{\partial p}\frac{\partial x_i}{\partial q}
\left(\frac{\partial u}{\partial p}\frac{\partial v}{\partial q}+
\frac{\partial u}{\partial q}\frac{\partial v}{\partial p}\right)
\right)
\\
&=
\frac1C \sum_i \PB{x_i}{u}\PB{x_i}{v}
\end{align*}
Here $C=\det(g)$ when written as a $2\times2$ matrix.
\begin{align*}
C = \sum_{ij}\left(
\frac{\partial x_i}{\partial p}
\frac{\partial x_i}{\partial p}
\frac{\partial x_j}{\partial q}
\frac{\partial x_j}{\partial q}
-
\frac{\partial x_i}{\partial p}
\frac{\partial x_j}{\partial p}
\frac{\partial x_i}{\partial q}
\frac{\partial x_j}{\partial q}
\right)
\end{align*}
which gives (\ref{Geo_C})
\end{proof}

Let $\man\subset\Real^n$ be a closed genus 0 symplectic surface, and
let $\Psi_\star\colon \man\mapsto{S^2}$ be a bijective symplectic
analytic diffeomorphism, and $\Psi^\star\colon \CoSS\mapsto\CoM$ be
the corresponding pullback map.  Let $(\theta,\phi)$ be the spherical
coordinates on $S^2$, then
$(p=\Psi^\star(\cos\theta),q=\Psi^\star(\phi))$ are conjugate
coordinates on $\man$. However these coordinates are not defined for
the whole of $\man$. More importantly, the noncommutative analogue of
$(p,q)$ do not have matrix representation. We can avoid this problem
by setting $J_0=\Psi^\star(\cos\theta)$,
$J_1=\Psi^\star(\sin\theta\cos\phi)$, and
$J_2=\Psi^\star(\sin\theta\sin\Phi)$.

The conformal factor $C$ in (\ref{Geo_C}) can now be written
\begin{align}
C &= \frac1{(1-J_0^2)}
\sum_{ij} \PB{x_j}{x_i}\PB{J_0}{x_i}
\left( J_1\PB{J_2}{x_j} - J_2\PB{J_1}{x_j} \right)
\label{Geo_CJ}
\end{align}

The two above expression are examples of the following the theorem:

\begin{theorem}
\label{thm_Geo}
Let $\man\subset\Real^n$ be a symplectic surface, and $u\colon
\man\mapsto\Real$ be a function that is derived from the metric on
$\man$ and its embedding, using only differentiation. Then we can find
an expression for $u$ using only the Poisson bracket, the embedding
functions $\Set{x_1,\ldots,x_n}$ and the conjugate coordinates
$\Set{p,q}$. If $\man$ is topologically the sphere then we can replace
$\Set{p,q}$ with $\Set{J_0,J_1,J_2}$.
\end{theorem}

\begin{proof}
Take the expression for $u$ and replace the metric with (\ref{Geo_C})
or (\ref{Geo_CJ}), 
and replace the derivatives using 
\begin{align*}
\frac{\partial u}{\partial q}=\PB{p}{u}
\,,\
\frac{\partial u}{\partial p}=-\PB{q}{u}
\end{align*}
or
\begin{align*}
\frac{\partial u}{\partial q}=\PB{J_0}{u}\,,\
\frac{\partial u}{\partial p}=(1-J_0^2)^{-1}
(J_1\PB{J_2}{u}-J_2\PB{J_1}{u})
\end{align*}
\end{proof}
Examples of such functions include the curvature and Laplacian, which
depend only on the metric, and the first and second fundamental forms,
which depend on the metric and the embedding.

\section{Examples}
\label{ch_eg}

\subsection{Heisenberg ANCG $\Heis_{2r,s}$} 
\label{ch_egHeis}

In section \ref{ch_Heis} we gave the details of the Heisenberg ANCG
$\Heis_{2r,s}$, including the Wick and normal orderings and their
corresponding star products. The Heisenberg algebra may be interpreted
as the noncommutative Euclidean flat space $\Real^{2r+s}$.  Clearly
the Heisenberg ANCG is its own coordinate chart. In section
\ref{ch_QG} we gave the quantum group based on $\Heis_{2r,s}$.

Because of the equation $[\elAp_i,\elAq_i]=i\varepsilon\delta_{ij}$, 
there do not exist any matrix representations of
$\Heis_{2r,s}$. There do however exist many (topologically
inequivalent) representations of $\Heis_{2r,s}$. 


\subsection{A Phase space $\man=T^\star S^2$}
\label{ch_egTS2}

Non-relativistic quantum mechanics is obtained via the
``quantisation'' of phase space. In our language this means finding an
ANCG $\Alg$ such that the corresponding manifold $\man=T^\star Q$ for
some configuration space $Q$, and such that the inherited Poisson
structure, is the canonical symplectic structure.

We give here an example corresponding to a free particle on a sphere,
so that $Q=S^2$ and $\man=T^\star S^2$. Note that, in order to keep
$\Alg$ algebraic, we require that we embed $T^\star S^2$ in
$\Real^8$, via the following embedding:
\begin{equation}
\begin{aligned}
\elOx_1 &= \sin\theta\,\cos\phi\,,\ &
\elOx_2 &= \sin\theta\,\sin\phi\,,\ &
\elOx_3 &= \cos\theta \,,\
\\
\elOx_4 &= \cos\theta\,\cos\phi \,,\ &
\elOx_5 &= \cos\theta\,\sin\phi \,,\ &
\elOx_6 &= \sin\theta\,,\ 
\\
\elOx_7 &= p_\theta  \,,\ &
\elOx_8 &= p_\phi  \,,\ &&
\end{aligned}
\label{egTS2_def_x}
\end{equation}
where $(\theta,\phi,p_\theta,p_\phi)$ is a coordinate chart for
$\man$, $(\theta,\phi)$ are the standard spherical coordinates, and
$p_\theta$ and $p_\phi$ there respective conjugate coordinates.

It is easy to show that the ANCG equivalent to the Heisenberg
quantisation of $\man=T^\star S^2$ is generated by
$\Set{\varepsilon,\elAx_1,\ldots,\elAx_8}$ with
$\pi(\elAx_i)=\elOx_i$.  The commutation relations are
\begin{equation}
\begin{array}{llll}{}
{[}\elAx_7,\elAx_1{]} = i\varepsilon\elAx_4 {\,,\ }&
{[}\elAx_7,\elAx_4{]} = -i\varepsilon\elAx_1 {\,,\ }&
{[}\elAx_8,\elAx_1{]} = -i\varepsilon\elAx_2 {\,,\ }&
{[}\elAx_8,\elAx_4{]} = i\varepsilon\elAx_5 {\,,\ }
\\{}
{[}\elAx_7,\elAx_2{]} = i\varepsilon\elAx_5  {\,,\ }&
{[}\elAx_7,\elAx_5{]} = -i\varepsilon\elAx_2 {\,,\ }&
{[}\elAx_8,\elAx_2{]} = i\varepsilon\elAx_1 {\,,\ }&
{[}\elAx_8,\elAx_5{]} = -i\varepsilon\elAx_4 {\,,\ }
\\{}
{[}\elAx_7,\elAx_3{]} = -i\varepsilon\elAx_6 {\,,\ }&
{[}\elAx_7,\elAx_6{]} = i\varepsilon\elAx_3 {\,,\ }&
{[}\elAx_8,\elAx_3{]} = 0 {\,,\ }&
{[}\elAx_8,\elAx_6{]} = 0 {\,,\ }
\\
\multicolumn{4}{l}{[\elAx_i,\elAx_j]=0\qquad\text{otherwise}} 
\end{array}
\label{egTS2_comm}
\end{equation}
and the immersion relations are
\begin{equation}
\begin{aligned}
\elAx_1^2 + \elAx_2^2 + \elAx_3^2 = 1\,,\qquad
\elAx_1^2 + \elAx_2^2 = \elAx_6^2 \,,
\\
\elAx_1\elAx_3 = \elAx_4\elAx_6 \,,\qquad 
\elAx_2\elAx_3 = \elAx_5\elAx_6 
\end{aligned}
\label{egTS2_imm}
\end{equation}
from which all the other immersion equations can be derived.

\vspace{1 em}

To give the standard Hamiltonian for a free particle on a sphere it
is necessary to enlarge $\Alg$ to include the generator $\elAx_9$ so
that $\elOx_9=1/\sin\theta$. Thus we must include the 
commutation relations
$[\elAx_7,\elAx_9]=-i\varepsilon\elAx_3\elAx_9^2$,
$[\elAx_9,\elAx_i]=0$ for $i\ne7$, and the immersion relation
$\elAx_6\elAx_9=1$. Topologically this is the noncommutative version
of the space $\man=T^\star(S^2\backslash\Set{N,S})$ where $N,S$ are
the two poles. The Hamiltonian is given by
$\elA{H}=\tfrac12 \elAx_7^2+\tfrac12\elAx_8^2\elAx_9^2$. Of
course this Hamiltonian is not unique and we can add any constant or
any multiple of $\varepsilon$ without effecting the classical dynamics. 

The Schroedinger representation is given by
$\Hil=L^2(S^2\backslash\Set{N,S})$, together with the inner product
$\innerprod{f,g}=\int_{S^2}\cnj{f}g\sin\theta\,d\theta\,d\phi$.  The
unbounded operators are given by
$\varphi(\elAx_i)f=\elOx_i f$, for $i=1,\ldots,6$ and
$\varphi(\elAx_7)f = i\varepsilon_\infty  \partial_\theta f$ and
$\varphi(\elAx_8)f = i\varepsilon_\infty  \partial_\phi f$.

A Heisenberg coordinate chart is given by
\begin{align*}
\Psi\colon\Alg\mapsto
\Heis_4
\left(
\begin{pmatrix}
{}^{-\infty}_{-\infty}
\\
{}^{-\infty}_{-\infty}
\end{pmatrix}
,
\begin{pmatrix}
{}^{\infty}_{\infty}
\\
{}^{\infty}_{\infty}
\end{pmatrix}
\right)
\qquad
\text{with coordinates}
\quad
(\elA{\theta},\elA{\phi},\elA{p_\theta},\elA{p_\phi})
\text{ where }
[\elA{p_\theta},\elA{\theta}] = [\elA{p_\phi},\elA{\phi}] =
i\varepsilon
\end{align*}
The immersion elements $\Set{ \elAx_1,\ldots,\elAx_8}$ are given by
(\ref{egTS2_def_x}), but replacing the unbolded with the bolded
symbols.


\subsection{Torus or Manin plane}
\label{ch_egT2}

For historical reasons the noncommutative torus is often called the
Manin plane or Weyl algebra. To place it in our language, $\Alg_{T^2}$
is generated by $\Set{\varepsilon,\elAx_1,\elAx_2,\elAx_3,\elAx_4}$
with $\elAx_1=\tfrac12(\elAu+\elAu^{-1})$,
$\elAx_2=\tfrac1{2i}(\elAu-\elAu^{-1})$,
$\elAx_3=\tfrac12(\elAv+\elAv^{-1})$,
$\elAx_4=\tfrac1{2i}(\elAv-\elAv^{-1})$. The relations are given by
\begin{align}
&\elAu^r\elAv^s=e^{i\varepsilon r s} \elAv^r\elAv^s\,, &
&\elAu^{-1}\elAu=\elAu\elAu^{-1}=\elAv^{-1}\elAv=\elAv\elAv^{-1}=1\,, &
&\elAu^\dagger=\elAu^{-1}\,, &
&\elAv^\dagger=\elAv^{-1} 
\end{align}
where $r,s=\pm1$. We can show that the first equation above is true
for all $r,s\in\Intg$.

There is a coordinate systems for the noncommutative torus given by
\begin{align}
\Psi_{T^2}\colon \Alg\mapsto\Heis_2
\left(
\begin{pmatrix}
{}^{-\infty}_{-\infty}
\end{pmatrix}
,
\begin{pmatrix}
{}^{\infty}_{\infty}
\end{pmatrix}
\right) \,;
&&
\Psi_{T^2}(\elAu)=e^{i\varepsilon\elAp} \,,&&
\Psi_{T^2}(\elAv)=e^{i\varepsilon\elAq} 
\end{align}

To get the Vey product we must use central ordering (theorem
\ref{th_Vey_cent} below) given by
$\Omega_V(\elOu^r\elOv^s)=\elAu^r\elAv^s e^{-i r s \varepsilon/2}$.

A normal ordering is given by
$\Omega_N(\elOu^r\elOv^s)=\elAu^r\elAv^s$. This produces the following
star product
\begin{align}
f\star_N g &=
\exp\left(-i\varepsilon 
\frac{\partial_2 }{\partial u}
\frac{\partial_1 }{\partial v}
\right)(f,g)
\end{align}

There is a matrix representation of $\Alg_{T^2}$ given with respect to
the basis $\Set{\rvec{0},\ldots,\rvec{N-1}}$
\begin{align}
\varphi_N(\elAu)\rvec{r}=e^{i r\varepsilon_N} \rvec{r}\,, \qquad
\varphi_N(\elAv)\rvec{0}=\rvec{N-1} \,,\
\varphi_N(\elAv)\rvec{r}=\rvec{r-1} \,,\ r=1,\ldots,N-1 
\end{align}
where $\varepsilon_N=1/N$. The trace map is therefore given by
\begin{align}
\Trn{N}(\elAu^r\elAv^s)=
\delta(r\,\text{mod}\,N)\,\delta(s\,\text{mod}\,N)
\end{align}
Therefore $\lim_{N\to\infty}\Trn{N}(\elA{y})$ exists for all
$\elA{y}\in\Alg_{T^2}$, and theorem \ref{thm_int} applies. However
$\Trn{N}$ is not analytic.

There is a Quantum Group structure for $\Alg_{T^2}$, suggested by
section \ref{ch_QG}, given by
\begin{align}
\begin{array}{llll}
\Delta(1)=1\otimes 1\,,\qquad 
&
\Delta(e^{i\varepsilon})=e^{i\varepsilon}\otimes
 e^{i\varepsilon},\qquad
&
\Delta(\elAu^r)=\elAu^r\otimes\elAu^r\,,\qquad 
&
\Delta(\elAv^r)=\elAv^r\otimes\elAv^r
\\
\epsilon(1)=1\,,&
\epsilon(e^{i\varepsilon})=1\,,&
\epsilon(\elAu^r)=1\,,&
\epsilon(\elAv^r)=1\,,
\\
S(1)=1\,,&
S(e^{i\varepsilon})=e^{-i\varepsilon}\,,&
S(\elAu^r)=\elAu^{-r}\,,&
S(\elAv^r)=\elAv^{-r}\,,
\end{array}
\end{align}
for all $r\in\Intg$.


\subsection{Surfaces of Rotation}
\label{ch_egSR}

These were first introduced in \cite{JG-SR} then expanded in
\cite{JG-TSR}. The ANCG, $\Algp$, are generated by
$\Set{\varepsilon,\elAx_1,\elAx_2,\elAx_3}$ and defined with respect
to a polynomial function $\rho\colon\Real^2\mapsto\Real$. The quotient
relations are given by
\begin{equation}
\begin{array}{c}
{[}\elAX_0,\elAX_+]=\varepsilon\elAX_+ \,, \qquad
{[}\elAX_0,\elAX_-]=-\varepsilon\elAX_- \,, \qquad
{[}\elAX_+,\elAX_-]=
\rho(\elAX_0-\varepsilon/2,\varepsilon)
-\rho(\elAX_0+\varepsilon/2,\varepsilon) \,,
\\
\elAX_+\elAX_- + \elAX_-\elAX_+ = 
\rho(\elAX_0-\varepsilon/2,\varepsilon)
+\rho(\elAX_0+\varepsilon/2,\varepsilon)
\end{array}
\label{egSR_qe}
\end{equation}
where $\elAx_1=\tfrac12(\elAX_+ + \elAX_-)$, $x_2=\tfrac1{2i}(\elAX_+ -
\elAX_-)$, $\elAx_3=\elAX_0$.

The topology of the corresponding $\man$ depends on the shape of the
curve $y(z)=\rho(z,0)$. If we let
$I_\rho(0)=\Set{z\in\Real|\rho(z,0)\ge0}$ then $I_\rho(0)$ is the
union of intervals. Assuming that $\rho(z,0)\ne0$ on the interior of
$I_\rho(0)$ then each bounded interval in $I_\rho(0)$ corresponds to a
disjoint submanifold topologically equivalent to the sphere. If one of
the intervals is either $\Set{z|-\infty<z<z_{\hi}}$ or
$\Set{z|z_{\lo}<z<\infty}$ then the corresponding submanifold is
topologically the disc. Finally if $I_\rho(0)=\Real$ then $\man$ is
topologically a cylinder.  For a $\rho$ with several maxima there may
be several intervals in $I_\rho(0)$, and therefore $\man$ is
disconnected. Replacing $\rho(z,\varepsilon)\to\rho(z,\varepsilon)+C$
may change the topology to $\man$. This is analysed in \cite{JG-TSR}.

If $\rho(z,\varepsilon)=z^2$ then $\man$ is not a manifold, but an
algebraic variety. However much of the analysis is still valid in this
case.

For $\varepsilon_0\in\Real$, $\varepsilon_0\ge0$ let
$I_\rho(\varepsilon_0)=\Set{z\in\Real|\rho(z,\varepsilon_0)>0}$. If
there exists $\varepsilon_N>0$ such that $I_\rho(\varepsilon_N)$ is a
bounded interval given by $I_\rho(\varepsilon_N)=
\Set{z|z_\lo(\varepsilon_N)<z<z_{\hi}(\varepsilon_N)}$ where
$N\varepsilon_N=z_\hi(\varepsilon_N)-z_\lo(\varepsilon_N)$, then there
is a $M_N(\Cmpx)$ representation of $\Algp$ given by
\begin{equation}
\begin{aligned}
\varphi_N(\elAX_0)\rvec{r}&=
\left(z_{\lo}(\varepsilon_N)+(r+\tfrac12)\varepsilon_N \right) \rvec{r} 
\\
\varphi_N(\elAX_+)\rvec{r}&=
\rho(z_{\lo}(\varepsilon_N)+(r+1)\varepsilon_N,\varepsilon_N)^\scrhalf
\rvec{r+1} 
\\
\varphi_N(\elAX_-)\rvec{r}&=
\rho(z_{\lo}(\varepsilon_N)+\varepsilon_N,\varepsilon_N)^\scrhalf 
\rvec{r-1} 
\\
\end{aligned}
\label{egSR_rep}
\end{equation}
If for some $\varepsilon_\infty>0$, $I_\rho(\varepsilon_\infty)$ is an
unbounded interval then there are infinite dimensional representations
of $\Algp$. It is easy to see that if $I_\rho(\varepsilon_0)$ is a
bounded interval for all $\varepsilon_0>0$ then the trace map is
defined, and furthermore it is analytic.

If $I_\rho(\varepsilon_0)$ is a single interval, possibly unbounded,
for all $\varepsilon_0>0$, then the coordinate system is given by
\begin{equation}
\begin{array}{c@{\qquad}c}
\Psi\colon \Alg_{\rho}\mapsto\Heis_2
\left(
\begin{pmatrix}
{}^{z_{\lo}}_{-\infty}
\end{pmatrix}
,
\begin{pmatrix}
{}^{z_{\hi}}_{\infty}
\end{pmatrix}
\right) \,;
&
\Psi(\elAX_0)=\elAp \,, 
\\
\Psi(\elAX_+)=
e^{i\elAq} (\rho(\elAp+\tfrac12\varepsilon,\varepsilon))^\scrhalf  \,,
&
\Psi(\elAX_-)=
e^{-i\elAq} (\rho(\elAp-\tfrac12\varepsilon,\varepsilon))^\scrhalf 
\end{array}
\label{egSR_coord}
\end{equation}
If $I_\rho(\varepsilon_\infty )$ is an unbounded interval we can
replace $z_{\lo}$ with $-\infty$ or $z_{\hi}$ with $+\infty$ or both.

\vspace{1 em}

Let $U_\rho=\Set{(z,\varepsilon_0)\in\Real^2\,|\,
\rho(z,\varepsilon_0)>0}$.  It is useful to enlarge $\Alg$ to the set
\begin{align}
\Algp= \Set{ \sum_{r=0}^\finite 
\elAX_+ f_r(\elAX_0,\varepsilon) + 
\sum_{r=0}^\finite
\elAX_- f_{-r}(\elAX_0,\varepsilon) }
\end{align}
where $f_r\colon U\mapsto\Cmpx$ is $C^\omega$ on the interior of
$U_\rho$. From (\ref{egSR_qe}) we have that
\begin{align}
f(\elAX_0,\varepsilon)\elAX_\pm=\elAX_\pm
f(\elAX_0\pm\varepsilon,\varepsilon). 
\end{align}
Because of this extension, we can talk about $\Algp$ even when
$\rho\colon U_\rho\mapsto\Real$ is bounded and $C^\omega$
on the interior of $U_\rho$.

As well as the homeomorphism giving the coordinate system, there are
isomorphisms between certain topologically equivalent noncommutative
surfaces of rotation. For example let $\Alg_{\rho_1}$ and
$\Alg_{\rho_2}$ be noncommutative surfaces of rotation with generators
$\varepsilon_1,\elAX_0,\elAX_+,\elAX_-$ and
$\varepsilon_2,\elAY_0,\elAY_+,\elAY_-$ respectively, such that
$\rho_1$, $\rho_2$ independent to $\varepsilon$, and both
$I_{\rho_1}=\Set{z|z^1_{\lo}<z<z^1_{\hi}}$ and
$I_{\rho_2}=\Set{z|z^2_{\lo}<z<z^2_{\hi}}$ are bounded then
\begin{align}
&
\Psi\colon \Alg_{\rho_1}\mapsto\Alg_{\rho_2}\,; 
\qquad
\Psi(\varepsilon_1)=\varepsilon_2 \,,
\qquad
\Psi(\elAX_0)=K( \elAY_0 - z^2_\lo)+ z^1_\lo \,,
\notag
\\
&
\displaystyle{
\Psi(\elAX_+)=\elAY_+ 
\left(
\frac{\rho_1(K(\elAY_0-z^2_\lo)+ z^1_\lo + \tfrac12\varepsilon_2)}
{\rho_2(\elAY_0+ \tfrac12\varepsilon_2)}
\right)^\scrhalf
} \,,
\qquad
\Psi(\elAX_-)=\Psi(\elAX_+)^\dagger
\label{egSR_Homs}
\end{align}
where 
\begin{align*}
K=\frac{(z^1_\hi-z^1_\lo)}{(z^2_\hi-z^2_\lo)}
\end{align*}

\vspace{1 em}

One possible ordering is the normal ordering given by 
\begin{align}
\Omega_N(\elOX_\pm^r f(\elOX_0)) &= \elAX_\pm^r f(\elAX_0)
\label{egSR_Om_N}
\end{align}
This ordering does not correspond to a differential star product with
$C_1:\SAlg\times\SAlg\mapsto\SAlg$ a first order operator. To see this
we note that
\begin{align*}
C_1(\elOX_+^r,\elOX_-^r)
&=
r\rho_\varepsilon(\elOX_0,0)\rho(\elOX_0,0)^{r-1}
+2^{-2r}(2r)!(n!)^{-1}
\rho_p(\elOX_0,0)\rho(\elOX_0,0)^{r-1}
\end{align*}
where $\rho_p$ and $\rho_\varepsilon$ are the partial differentiation
of $\rho(\elO{p},\varepsilon)$ with respect to the the first and
second arguments respectively.

In order to get the Vey product we need the \defn{central ordering}
which is defined with respect to the coordinates $(\elAp,\elAq)$. If
$\Psi$ is the coordinate homomorphism (\ref{egSR_coord}) then
\begin{align}
\Omega_C(\elOu)=\elAu \,,\text{ where }
\Psi^0(\elOu)=e^{i r\elOq}f(\elOp)
\text{ and }
\Psi(\elAu)=e^{i r\elAq}f(\elAp+r\varepsilon/2)
\label{egSR_Om_W}
\end{align}
In terms of the elements of $\Algp$ we can show that
\begin{equation}
\begin{aligned}
\Omega_C(\elOX_+^r f(\elOX_0)) &=
\elAX_+^r 
\left(
\frac
{(\rho(X_0+\tfrac12r\varepsilon,\varepsilon))^{r-1}}
{
\rho(\elAX_0+\tfrac32\varepsilon,\varepsilon)
\rho(\elAX_0+\tfrac52\varepsilon,\varepsilon)\cdots
\rho(\elAX_0+\tfrac{2r-1}2\varepsilon,\varepsilon)
}
\right)^{1/2} f(\elAX_0+\tfrac12r\varepsilon)
\\
\Omega_C(\elOX_-^r f(\elOX_0)) &= 
\elAX_-^r 
\left(
\frac
{(\rho(X_0-\tfrac12r\varepsilon,\varepsilon))^{r-1}}
{
\rho(\elAX_0-\tfrac32\varepsilon,\varepsilon)
\rho(\elAX_0-\tfrac52\varepsilon,\varepsilon)\cdots
\rho(\elAX_0-\tfrac{2r-1}2\varepsilon,\varepsilon)
}
\right)^{1/2} f(\elAX_0-\tfrac12r\varepsilon)
\end{aligned}
\label{egSR_Om_C}
\end{equation}

\begin{theorem}
\label{th_Vey_cent} 
The central ordering is compatible with the Wick ordering under the
Heisenberg coordinate homeomorphism:
\begin{align}
\Psi\circ\Omega_C &= \Omega_W\circ\Psi^0
\end{align}
where $\Psi$ is given by (\ref{egSR_coord}) and $\Omega_W$ by
(\ref{Heis_Om_W}).  Hence the central ordering gives the Vey product.
\end{theorem}
\begin{proof}
From (\ref{Heis_N2S}) we have
\begin{align*}
\Omega_W(e^{ib\elOq}\elOp^a) 
&=
\sum_{s=0}^\infty \elA{S}(a,s)\frac{(bi)^s}{s!}
=
\sum_{s=0}^\infty
\sum_{r=0}^a
\frac{b^s i^{s-r}(\tfrac{-\varepsilon}2)^r a!}
{(a-r)!r!(s-r)!}
\elAp^{a-r}\elAq^{s-r}
\\
&=
\sum_{r=0}^\infty
\sum_{t=0}^a
\frac{(bi)^t}{t!}
\frac{b^r (\tfrac{-\varepsilon}2)^r a!}
{(a-r)!r!}
\elAp^{a-r}\elAq^{t}
=
(\elAp-b\varepsilon/2)^a
e^{i\elAq}
\end{align*}
Hence (\ref{egSR_Om_W}). Using theorem  \ref{th_Hom_Ord_star} shows
that the star product must be Vey.

We can also prove that the central ordering gives the
Vey product directly.  Let $\elOu=e^{inq}f(p)$ and $\elOv=e^{imq}g(p)$
then from the definition of the Vey product we have
\begin{align*}
\elOu\star\elOv &=
\sum_{s=0}^\infty \frac{(i\varepsilon/2)^s}{s!}P^s(F,G) 
\\
&=
\sum_{s=0}^\infty \frac{(i\varepsilon/2)^s}{s!}
\left(\frac{\partial_1}{\partial p}\frac{\partial_2}{\partial q}
-\frac{\partial_1}{\partial p}\frac{\partial_2}{\partial q}\right)^s 
(u,v)
\\
&=
\sum_{t=0}^\infty \frac{(i\varepsilon/2)^t}{t!}
\left(\frac{\partial_1}{\partial p}\frac{\partial_2}{\partial
q}\right)^t 
\sum_{r=0}^\infty \frac{(i\varepsilon/2)^r}{r!}
\left(\frac{\partial_2}{\partial p}\frac{\partial_1}{\partial
q}\right)^r 
(u,v)
\\
&=
e^{iq(n+m)}f(p+m\varepsilon/2)g(q-n\varepsilon/2)
\end{align*}
where $P(u,v)=\PB{u}{v}$ and $(\partial_1/\partial p)$ refers to
differentiating with respect to $u$ and $(\partial_2/\partial p)$
refers to differentiating with respect to $v$.  Thus
\begin{align*}
\Omegas_C(u\star v)&=
e^{i\elAq(n+m)}f(\elAp+m\varepsilon+n\varepsilon/2)g(q+m\varepsilon/2)
\\
&=e^{i\elAq n}f(\elAp+n\varepsilon/2)
e^{i\elAq m}g(q+m\varepsilon/2)
\\
&=
\Omegas_C(u)\Omegas_C(v)
\end{align*}
\end{proof}

In general we can only give the quantum group structure in terms of a
formal expansion.  By directly applying (\ref{QG_Heis}) on
(\ref{egSR_coord}), we get
\begin{align}
&
\Delta(\varepsilon)=0 \,,\qquad\qquad
\Delta(\elAX_0)=1\otimes \elAX_0 + \elAX_0\otimes 1 \,,
\label{egSR_QG}
\\
&
\Delta(X_+) = \sum_{a,b=0}^{\infty}\sum_{s=0}^a\sum_{t=0}^b
\frac{\alpha_{ab} a!b!}{s!(a-s)!t!(b-t)!}
\varepsilon^s
\elAX_+ \rho(\elAX_0\!+\!\tfrac12\varepsilon,\varepsilon)^{-\scrhalf} 
\elAX_0^t
\otimes 
\varepsilon^{a-s}
\elAX_+ \rho(\elAX_0\!+\!\tfrac12\varepsilon,\varepsilon)^{-\scrhalf} 
\elAX_0^{b-t} \,,
\notag
\\
&
\Delta(X_-) = \sum_{a,b=0}^{\infty}\sum_{s=0}^a\sum_{t=0}^b
\frac{\alpha_{ab} a!b!}{s!(a-s)!t!(b-t)!}
\varepsilon^s
\rho(\elAX_0\!+\!\tfrac12\varepsilon,\varepsilon)^{-\scrhalf} 
\elAX_0^t
\elAX_- 
\otimes 
\varepsilon^{a-s}
\rho(\elAX_0\!+\!\tfrac12\varepsilon,\varepsilon)^{-\scrhalf} 
\elAX_0^{b-t}
\elAX_- \,,
\notag
\end{align}
\begin{align*}
&
\epsilon(\varepsilon)=0 \,,\qquad 
\epsilon(\elAX_0)=0 \,,\qquad 
\epsilon(\elAX_+)=\rho(0,0) \,,\qquad
\epsilon(\elAX_-)=\rho(0,0)\,,
\\
&
S(\varepsilon)=-\varepsilon \,,\qquad
S(\elAX_0)=-\elAX_0 \,,
\\
&
S(\elAX_+)=\elAX_- 
\rho(-\elAX_0+\tfrac12\varepsilon,-\varepsilon)^\scrhalf
\rho(\elAX_0-\tfrac12\varepsilon,\varepsilon)^{-\scrhalf} \,,
\\
&
S(\elAX_-)=
\rho(-\elAX_0+\tfrac12\varepsilon,-\varepsilon)^\scrhalf
\rho(\elAX_0-\tfrac12\varepsilon,\varepsilon)^{-\scrhalf}
\elAX_+ \,,
\end{align*}
where the Taylor expansions of $\rho^\scrhalf$ is given by
\begin{align*}
\rho(u+\tfrac12\varepsilon,\varepsilon)^\scrhalf=\sum_{a,b=0}^\infty
\alpha_{ab}\varepsilon^a u^b
\end{align*}
Note these simplifies a little if
$\rho(u,\varepsilon)=\rho(-u,-\varepsilon)$. Clearly for the image of
$\Delta$ in (\ref{egSR_QG}) to be a polynomial requires that
$\rho(u+\tfrac12\varepsilon,\varepsilon)^\scrhalf$ is a
polynomial. This implies that $\man_\rho$ is topologically the
cylinder. (Thus excluding the sphere.)  Examples of such $\rho$
include $\rho(u,\varepsilon)=1$ and
$\rho(u,\varepsilon)=(u^2+1)^2$,


\subsection{The Sphere}
\label{ch_egS2}

The noncommutative sphere has been studied by many
authors \cite{Madore_bk,JG-TS2,JG-S2,Cahen1}. It is an example of a
noncommutative surface of rotation with
\begin{align}
\rho(z,\varepsilon)=R^2-z^2+\varepsilon^2/4
\end{align}
where $R\in\Real$ gives the radius of the embedded sphere.  By looking
at the commutation relations part of (\ref{egSR_qe}) we see that
$\Set{\elAx_1,\elAx_2,\elAx_3}$ obey the commutation relations of the
Lie algebra $su(2)$, given by
$[\elAx_i,\elAx_j]=i\varepsilon\,\epsilon_{ijk} \elAx_k$. The immersion
equation of (\ref{egSR_qe}) gives the Casimir
$\elAx_1^2+\elAx_2^2+\elAx_3^2=R^2$. As a result \cite{JG-S2} both
$\Alg_{S^2}$ and $\CAlg_{S^2}$ are infinite dimensional
representations of $su(2)$, with $\CAlg_{S^2}$ being the enveloping
algebra.

All the results for noncommutative surface of rotation now carry over,
including the Vey or central ordering.
The finite dimensional unitary representation of $\Alg_{S^2}$ given by
(\ref{egSR_rep}) reduces to
the standard representations of $su(2)$ 
\begin{equation}
\begin{aligned}
\varphi_N(\varepsilon) &= \varepsilon_N=2R(N^2-1)^{-\scrhalf}
\\
\varphi_N(\elAX_0) \rvec{m} 
&= \varepsilon_N (m-\tfrac{N-1}2) \rvec{m} 
\\
\varphi_N(\elAX_+) \rvec{m} 
&= 
\varepsilon_N (N-m-1)^{\scrhalf} (m+1)^{\scrhalf} \rvec{m+1} 
\\
\varphi_N(\elAX_-) \rvec{m} 
&= 
\varepsilon_N (N-m)^{\scrhalf} (m)^{\scrhalf} \rvec{m-1} 
\end{aligned}
\label{egS2_rep}
\end{equation}
where $N\in\natnum$. 

\vspace{1 em}

Additionally the noncommutative sphere has a
specific basis $\elA{P}^n_m$ for $m,n\in\Intg$, $n\ge0$ and $|m|\le
n$, and specific ordering $\Omega_{S^2}$. This was given in detail in
\cite{JG-S2,JG-SR}, although it has a slightly different
normalisation there.  The elements $\elA{P}^n_m$ are defined via
\begin{align}
\elA{P}^m_n &=
\alpha^{m-n}\varepsilon^{m-n}
\left(\frac{(n+m)!}{(2n)!\,(n-m)!}\right)^{\scrhalf} 
(\ad{\elAX_-})^{n-m}(\elAX_+{}^n)
\label{egS2_def_Pmn}
\end{align}
where $\ad{\elAu}(\elAv)=[\elAu,\elAv]$.  When written as a formally
tracefree symmetric polynomial in $(\elAX_0,\elAX_+,\elAX_-)$,
$\elA{P}^m_n$ is a homogeneous polynomial of order $n$ and is
independent of $R$ and $\varepsilon$. This justifies
(\ref{egS2_OmS2}), as the spherical harmonics can also be written as
formally tracefree symmetric polynomials.

There is a sesquilinear form on $\Alg_{S^2}$ defined by
$\innerprod{\elAu,\elAv}=\pi_0(\elAu^\dagger\elAv)$ where
$\pi_0(\elAu)$ is the coefficient of $\elAu$ independent of $\elAx_i$
when $\elAu$ is written as a formally tracefree symmetric polynomial.
The sesquilinear form is related to the trace via
\begin{align}
\pi_0(\elAu)=\TrA(\elAu)
\end{align}
With respect to this sesquilinear form the basis  elements
$\elA{P}^n_m$ are orthogonal.
Each $\elA{P}^m_n$ is an eigenvector of the operators $\ad{\elAX_0}$ and
$\Delta=\ad{\elAX_0}^2 +
\tfrac12(\ad{\elAX_+}\ad{\elAX_-}+\ad{\elAX_-}\ad{\elAX_+})$:
\begin{align}
\ad{\elAX_0} \elA{P}^m_n &= \varepsilon m \elA{P}^m_n 
\label{egS2_ez_Pmn} 
\\
\Delta \elA{P}^m_n &= \varepsilon^2 n(n+1) \elA{P}^m_n 
\label{egS2_Del_Pmn} 
\end{align}
The ladder operators $\ad{\elAX_+},\ad{\elAX_-}$ increase or decrease
$m$:
\begin{align}
\ad{\elAX_\pm} \elA{P}^m_n &= \alpha\varepsilon (n\mp m)^{\scrhalf} 
(n\pm m+1)^{\scrhalf} \elA{P}^{m\pm 1}_n 
\label{egS2_epm_Pmn}
\end{align}
and the ``normal'' of $\elA{P}^m_n$ is given by
\begin{align}
\innerprod{{\elA{P}^m_n},{\elA{P}^m_n}}=
\norm{\elA{P}^m_n}^2 &=
\alpha^{2n}
\frac{(n!)^2}{(2n+1)!} 
\prod_{r=1}^n(4R^2+\varepsilon^2(1-r^2))
\label{egS2_Norm}
\end{align}
The product of two basis elements is given in terms of Wigner $6j$
symbols:
\begin{align}
\elA{P}^{m_1}_{n_1}
\elA{P}^{m_2}_{n_2}
&=
\sum_{n=|n_1-n_2|}^{n=n_1+n_2}
\CG{n_1}{n_2}{n}{m_1}{m_2}{m_1+m_2}
\RMP{n_1}{n_2}{n}
\elA{P}_{n}^{m_1+m_2} 
\label{egS2_RM_def}
\end{align}
where $\CG{n_1}{n_2}{n}{m_1}{m_2}{m_1+m_2}$ is the Clebsh-Gordon
coefficient, and the reduced matrix element $\RMP{n_1}{n_2}{n}$ is given
by
\begin{align}
\RMP{n_1}{n_2}{n}
&=
(-1)^{N+1+n_1+n_2}
\frac{\norm{\elA{P}^{m_1}_{n_1}}\norm{\elA{P}^{m_2}_{n_2}}}
{\norm{\elA{P}^{m_1+m_2}_{n}}}
(N)^\scrhalf 
(2n_1+1)^\scrhalf (2n_2+1)^\scrhalf 
\WsixJ{\tfrac{N-1}{2}}{n_1}{\tfrac{N-1}{2}}{n_2}{\tfrac{N-1}{2}}{n} 
\label{egS2_RM_res}
\end{align}
where $N=(4R^2\varepsilon^{-2}+1)^\scrhalf$ and the symbol in the
curly brackets is Wigner's 6-$j$ coefficient. Note the right hand side
of (\ref{egS2_RM_res}) is only defined when $N\in\natnum$. 

The image of $\varphi_N(\elA{P}^n_m)$ is a Wigner Operator. This must
be written in half integer notation, where $2k+1=N$ and
$j=-k,-k+1,\ldots,k$.
\begin{align}
\varphi_N(\elA{P}^m_n)\rvec{k,j}
&=
(-1)^{n} 
\norm{\elA{P}^m_n}
(2n+1)^\scrhalf 
\Wigner{2n}{n}{n+m}
\rvec{k,j}
\label{egS2_Wig_op}
\end{align}

As well as the normal ordering and the central ordering there is a
``Wick-like'' ordering $\Omega_{S^2}$ is given by
\begin{align}
\Omega_{S^2}(\psi^m_n)=
(-1)^n \frac{((2n+1)!)^\scrhalf}{n!(2R)^n}
\elA{P}^n_m
\label{egS2_OmS2}
\end{align}
Hence $\elA{P}^n_m$ may be thought of as the noncommutative analogue
of spherical harmonics.  We also note that this ordering is compatible
with the trace map since $\Trn{N}(\Omega_{S^2}(\elOu))$ is independent
of $N$.  A closed formula for the corresponding star product is being
searched.  It is known \cite{Cahen1} that it is not a differential
star product.


\subsection{Complex and Other Planes}
\label{ch_egC}

\begin{lemma}
Let $\Falg=\Falginf$ be generated by $\elFx,\elFy,\varepsilon$. Let
$\QCP\in\Falg$ be any element. Then there is a noncommutative plane
given by $\Alg=\Falg/\Ideal\Set{[\elFx,\elFy]\sim i\varepsilon\QCP}$,
if $\elA{c}=\elA{c}^\dagger$ where $\elA{c}=\q(\QCP)$.
\end{lemma}
\begin{proof}
Trivial.
\end{proof}

As usual we set $\elAx,\elAy\in\Alg$ to be the images of $\elFx,\elFy$
under the quotient.  We note that in general this procedure does not
produce an ANCG for $n\ge2$. This is because in general the Jacobi
identity is not satisfied. This is a rich source of ANCGs.  If we set
$\elAz=\elAx+i\elAy$ and $\cnj{\elAz}=\elAx+i\elAy$ then this is often
called the noncommutative complex plane. An example is given in
\cite{Klimek_Lesn1}.


\section{An Application: Finite Models of compact surfaces}
\label{ch_Mod}

In this section we give a finite element method for analysing surfaces
based on expansions in spherical harmonics.  As mentioned in the
introduction, this method is based on noncommutative geometry and
hence there is an error introduced depending on the order of
multiplication. However, the result is associative.

\vspace{1 em}

Let us assume that $\man$ is a surface of genus $0$ and we have the
diffeomorphism $\Psi_\star\colon S^2\mapsto\man$. From $\Psi_\star$ we
generate the pull back map $\Psi^\star\colon \CoM\mapsto\CoSS$.

To convert the functions $u\colon \man\mapsto\Cmpx$ into matrices we
would ideally use the homomorphism
\begin{align*}
\Phi_N\colon \CoM\mapsto M_N(\Cmpx)\,;\qquad
\Phi_N=\varphi_N\circ\Omega_{S^2}\circ\Psi^\star
\end{align*}
where $\Omega_{S^2}$ is given in (\ref{egS2_OmS2}), and $\varphi_N$ is
given in (\ref{egS2_rep}).  However, in general, $\Psi^\star(u)\in
\CoSS$ does not belong to $\OAlg_{S^2}$; that is, a finite sum
of spherical harmonics.  As mentioned in section \ref{ch_new} we can
still define the image of $\Phi_N$ via
\begin{align}
\Phi_N(u)=
\sum_{n,m}
\left(\int_{S^2}
\cnj{\psi^m_n} \Psi^\star(u) \sin\theta d\theta d\phi
\right)
(-1)^n \frac{((2n+1)!)^\scrhalf}{n!(2R)^n}
\varphi_N(\elA{P}^n_m)
\end{align}
We can calculate  $\varphi_N(\elA{P}^n_m)$ using (\ref{egS2_Wig_op}).
Since we have a loss of information when converting from functions to
matrices we cannot expect an inverse map. However the ``one sided
inverse'' to $\Phi_N$ is given by
\begin{align}
\Upsilon_N\colon M_N(\Cmpx)\mapsto\CoM\,;\qquad
\Upsilon_N(u_N)=
\sum_{n,m}
\frac{\Trn{N}(\varphi_N(\elA{P}^n_m)^\dagger u_N)}
{\varphi_N(\norm{\elA{P}^n_m}^2)}
(-1)^n \frac{n!(2R)^n}{((2n+1)!)^\scrhalf}
\Psi^{\star -1}(\psi^n_m)
\end{align}
It is easy to show that these satisfy 
\begin{align}
&\Phi_N\circ\Upsilon_N=1_{M_N(\Cmpx)} 
\\
&u-\Upsilon_N(\Phi_N(u))=
\sum_{n=N}^\infty\sum_{m=-n}^n
\left(\int_{S^2}
\cnj{\psi^m_n} \Psi^\star(u) \sin\theta d\theta d\phi
\right)
\end{align}
If $\Set{x_1,x_2,x_3}$ are the immersion coordinates of $\man$ then
$\Set{\Phi_N(x_1),\Phi_N(x_2),\Phi_N(x_3)}$ encode the geometry of
$\man$ into matrices. Other ``external'' information on $\man$; that
is any function $h\colon\man\mapsto\Cmpx$ (for example representing
density,) is also encoded as $\Phi_N(h)$.

\vspace{1 em}

The next step is to convert the expression for the desired result, in
terms of a matrix expression. For this we employ theorem
\ref{thm_Geo} and theorem \ref{thm_int}. Theorem \ref{thm_Geo} states
that any differentiation can be written in terms of the Poisson
bracket. We can therefore use (\ref{Poi_def_Poi}) to give the
differentiation in terms of a commutator. Theorem \ref{thm_int} states
that we can replace integration with the trace.  Combining these we
show that if the result can be expressed solely in terms of the known
functions via integration and differentiation then we can rewrite the
expression as matrix operations. Applying this we obtain the result as
a matrix. Finally we apply $\Upsilon_N$ to obtain an approximate
result. Clearly the rate of convergence for this algorithm depends on
the rate of convergence of the modular expansion of the functions
$\Set{x_1,x_2,x_3}$ and the external functions $h_s$.


\section{Discussion}
\label{ch_Disc}

We have given a consistent definition of an algebra $\Alg$ in terms of
noncommuting coordinates of an immersion space. When a parameter
$\varepsilon$ is set to zero, we obtain the commutative algebra
$\OAlg$ of functions on an algebraic manifold $\man$. This $\OAlg$ is
a subalgebra of $\CoM$, which is dense if $\man$ is compact. We have
shown that $\man$ inherits a Poisson structure as the limit of the
commutator. If we give $\Alg$ an ordering then we obtain a star
product on $\man$. We have define homomorphism and isomorphisms
between noncommutative geometries. By mapping one noncommutative
geometry to the Heisenberg algebra, we have given an analogue of the
coordinate chart and have given $\Alg$ a quantum group
structure. Noncommutative versions of $\Real^n$, $T^\star S^2$, $T^2$,
$S^2$ and surfaces of rotation have been developed. The metric has
been extended to noncommutative geometry and used to give
an application of noncommutative geometry to the numerical analysis
of surfaces.

\vspace{1 em}

One of the principle challenges is to enlarge $\Alg$ so that
$\OAlg=\CoM$. This would enable us to generalise theorem
\ref{th_Ord_star_gives_A} and say that equation
(\ref{intro_star_apprx}) is an exact equivalence. We have already
suggested how to partially enlarge $\Alg$ for some examples such as
the surfaces of rotation and flat space.  One possibility is to use an
ordering to define $\SAlg$. If this ordering is chosen so that (1) we
can extend the domain of $\Omega$ to $\CoM$ and (2) the star product
was a differential star product then we would have such an extension.
An alternative, would be an intrinsic definition of an algebraic
noncommutative geometry using the coordinate charts described in
section \ref{ch_Heis}. This would require a definition of analytic
continuation. It would make the definition of a noncommutative
manifold independent of the immersion and similar in spirit to the
definition of a standard manifold.

We note that if $\Alg$ has a Banach or $c^\star$ structure, then this
could be used to complete $\Alg$. However in general such a structure
does not exists.  Thus this approach differs from that of Alan Connes,
who investigated an alternative definition of a noncommutative
geometry $\Alg$ so that it was a $c^\star$ algebra. As a result he
sets the maps $\pi$ and $\Omega$ so that, \cite[page 156]{Connes2}
\begin{align}
&\lim_{\varepsilon\to0}
(\Omega(u)+\lambda\Omega(v)-\Omega(u+\lambda v))=0
\,,
\qquad \forall u,v\in\OAlg\,,\quad\lambda\in\Cmpx
\label{connes_1}
\\
&\lim_{\varepsilon\to0}
(\Omega(u)\Omega(v)-\Omega(u v))=0
\,,
\qquad \forall u,v\in\OAlg
\label{connes_2}
\\
&\lim_{\varepsilon\to0}
(\Omega(u^\dagger)-\Omega(u)^\dagger)=0
\,,
\qquad \forall u\in\OAlg
\label{connes_3}
\end{align}
All these are true for the definition of $\Omega$ in this article
since $\pi\circ\Omega=1_{\OAlg}$. However (\ref{connes_1}) is true for
all $\varepsilon$ not just in the limit, and (\ref{connes_3}) is true
for all $\varepsilon$ if $\Omega$ is a unitary ordering.

\vspace{1 em}

As mentioned in the introduction, there are many ways of defining the
analogue of a tangent vector field, and we would like to extend the
definition of a vector field given in \cite{JG-TS2,JG-TSR} for spheres
and surfaces of rotation, to that of a general algebraic
noncommutative geometry.

\vspace{1 em}

Considering some of the physical applications of this theory;
as mentioned, noncommutative geometry has been suggested as a
candidate for quantum gravity. Since the classical spacetime inherits
a Poisson structure from the noncommutativity of $\Alg$, we should
apply this procedure to spacetimes such as the Schwarzchild black hole
where there is a ``natural'' Poisson structure arising from the
Killing-Yano tensors. This will enable one to study the suggestion by
'tHooft and others that the event horizon should contain only a finite
quantity of information.

An alternative applications is given in \cite{JG-RWT-Q}, where Gratus
and Tucker use an algebra, based on the noncommutative surface of
rotation, to describe a Q-brane, a possible model for states of matter.
They also suggest how to interpret $\man$ as a phase space
even when $\man$ is not a cotangent bundle.

Finally we would like to demonstrate situations in the real world of
mathematical modelling, where the method outlined in section
\ref{ch_Mod}, is more efficient than standard approaches.

\subsection*{Acknowledgement}

The author would like to thank Robin Tucker and Marianne
Karlsen for their suggestions and help in the preparation of this
article, and the physics department of Lancaster University for their
facilities.


\end{document}